\documentclass[10pt]{article}
\newcommand{\heuteIst}{Oct 9, 2001 }
\usepackage{amsmath,amsthm,amsfonts,latexsym,amscd,amssymb}
\newcommand{\href}[2]{#2}

\scrollmode

\swapnumbers
\theoremstyle{plain}
\newtheorem{theorem}{Theorem}[section]
\newtheorem{lemma}[theorem]{Lemma}

\newtheorem{proposition}[theorem]{Proposition}

\theoremstyle{definition}
\newtheorem{definition}[theorem]{Definition}

\newtheorem{notation}[theorem]{Notation}

\theoremstyle{remark}
\newtheorem{remark}[theorem]{Remark}

\newcommand{\reals}{\mathbb{R}}

\newcommand{\complexs}{\mathbb{C}}

\newcommand{\naturals}{\mathbb{N}}

\newcommand{\integers}{\mathbb{Z}}
\newcommand{\rationals}{\mathbb{Q}}
\newcommand{\Q}{\mathbb{Q}}

\newcommand{\1}{\mathbf{1}}  
\DeclareMathOperator{\id}{id}
\newcommand{\boundary}[1]{\partial#1}

\newcommand{\abs}[1]{\left\lvert#1\right\rvert} 
\newcommand{\norm}[1]{\lVert#1\rVert}

\newcommand{\tensor}{\otimes}
\newcommand{\into}{\hookrightarrow}

\newcommand{\iso}{\cong}

\newcommand{\NeumannN}{\mathcal{N}}

\newcommand{\fundamentalDomain}{\mathcal{F}}
\newcommand{\innerprod}[1]{\langle #1 \rangle}
\DeclareMathOperator{\im}{im}      
\DeclareMathOperator{\vol}{vol}    
\DeclareMathOperator{\area}{area}
\DeclareMathOperator{\diam}{diam}  

\DeclareMathOperator{\cone}{cone}

\DeclareMathOperator{\tr}{tr}
\DeclareMathOperator{\diag}{diag}
\DeclareMathOperator{\Det}{det}

\DeclareMathOperator{\interior}{int}

\DeclareMathOperator{\sign}{sign}
\DeclareMathOperator{\Star}{star}  

\newcommand{\Kommentar}[1]{}
\newcommand{\forget}[1]{}

{\catcode`@=11\global\let\c@equation=\c@theorem}

\allowdisplaybreaks[2]

\newcommand{\RAgroups}{{\mathcal{G}}}
\DeclareMathOperator{\pr}{pr}
\newcommand{\squarematrix}[4]                
{                                            
\begin{pmatrix} #1 & #2 \\ #3 &
#4
\end{pmatrix}
}

\begin{document}

\date{Last edited \heuteIst or later --- last compiled: \today}

\title{Approximating $L^2$-signatures by their compact analogues}
\author{ Wolfgang
L\"uck\thanks{\noindent email: \protect\href{wolfgang.lueck@math.uni-muenster.de}{wolfgang.lueck@math.uni-muenster.de}\protect\\
www: ~\protect\href{http://www.math.uni-muenster.de/u/lueck/org/staff/lueck/}{http://www.math.uni-muenster.de/u/lueck/org/staff/lueck/}\protect}\\
Fachbereich Mathematik\\ Universit\"at M\"unster\\
Einsteinstr.~62\\         48149 M\"unster \and
Thomas Schick\thanks{\noindent e-mail:
\protect\href{thomas.schick@math.uni-muenster.de}{thomas.schick@math.uni-muenster.de}\protect\\
www:~\protect\href{http://www.math.uni-muenster.de/u/schickt/}{http://www.math.uni-muenster.de/u/schickt/}\protect\\
Research partially carried out during a stay at Penn State university
funded by the DAAD}\\ Fakult\"at f\"ur Mathematik\\ Universit\"at
G\"ottingen\\ Bunsenstrasse 3\\
37073 G\"ottingen\\ Germany}
\maketitle

\typeout{--------------------  Abstract  --------------------}

\begin{abstract}
  Let $\Gamma$ be a group together with a sequence
  of normal subgroups $\Gamma\supset \Gamma_1\supset \Gamma_2\supset\dots$
  of finite index $[\Gamma:\Gamma_k]$ such
  that $\bigcap_k \Gamma_k=\{1\}$.
  Let $(X,Y)$ be a (compact) $4n$-dimensional Poincar\'e
  pair and $p: (\overline{X},\overline{Y}) \to (X,Y)$ be a
  $\Gamma$-covering, i.e. normal
  covering with $\Gamma$ as deck transformation group.
  We get associated $\Gamma/\Gamma_k$-coverings $(X_k,Y_k) \to (X,Y)$.
  We prove that
  \begin{equation*}
    \sign^{(2)}(\overline{X},\overline{Y}) = \lim_{k\to\infty}
    \frac{\sign(X_k,Y_k)}{[\Gamma : \Gamma_k]},
  \end{equation*}
  where $\sign$ or $\sign^{(2)}$
  is the signature or $L^2$-signature, respectively, and the
  convergence of the right side for
  any such sequence $(\Gamma_k)_{k \ge 1}$ is part
  of the statement.


  If $\Gamma$ is amenable, we prove in a similar way an
  approximation theorem for $\sign^{(2)}(\overline{X},\overline Y)$ in
  terms of the
  signatures of a regular exhaustion of $\overline{X}$.
\end{abstract}

\noindent
Key words: $L^2$-signature, signature,
covering with residually finite deck transformation group, amenable exhaustion.
\\
2000 mathematics subject classification: 57P10, 
57N65, 
58G10 

\typeout{--------------------   Introduction  --------------------}

\setcounter{section}{-1}
\section{Introduction}
\label{sec:intro}

Throughout most of this paper we will use the following conventions.
We fix a group $\Gamma$, first together
with a sequence of normal subgroups
$\Gamma\supset \Gamma_1\supset \Gamma_2\supset\dots$
of finite index $[\Gamma:\Gamma_k]$ such that $\bigcap_k \Gamma_k=\{1\}$.
(Provided that $\Gamma$ is countable,
$\Gamma$ is residually finite if and only if such a
sequence $(\Gamma_k)_{k \ge 1}$ exists.) Moreover,
given a $\Gamma$-covering $p: \overline{X} \to X$,
i.e. a normal covering with $\Gamma$ as group of
deck transformations, we will denote the associated
$\Gamma/\Gamma_k$-coverings by $X_k:=\overline{X}/\Gamma_k \to X$
and for a subspace $Y \subset X$ let
$\overline{Y} \subset \overline{X}$
and $Y_k \subset X_k$ be the obvious pre-images.

One of the main results of the paper is
\begin{theorem}\label{res_approxi}
Let $(X,Y)$ be a  $4n$-dimensional Poincar\'e pair.
Then the sequence
$\left(\sign(X_k,Y_k)/[\Gamma : \Gamma_k]\right)_{k \ge 1}$
converges and
  \begin{equation*}
   \lim_{k\to \infty}    \frac{\sign(X_k,Y_k)}{[\Gamma : \Gamma_k]} =
   \sign^{(2)}(\overline{X},\overline{Y}).
  \end{equation*}
\end{theorem}

Some explanations are in order.
An $l$-dimensional \emph{Poincar\'e pair}
$(X,Y)$ is a pair of finite $CW$-complexes $(X,Y)$
with connected $X$ together with a
so called fundamental class $[X,Y] \in H_l(X;\Q)$
such that for the universal covering, and hence for any
$\Gamma$-covering $p: \overline{X} \to X$, the
Poincar\'e $\Q \Gamma$-chain map induced
by the cap product with (a representative of) the fundamental class
$$\cdot \cap [X,Y]: C^{l-*}(\overline{X},\overline{Y}) \to C_*(\overline{X})$$
is a $\Q \Gamma$-chain homotopy equivalence. Because we are working
with free finitely generated $\rationals\Gamma$-chain complexes, this
is the same as saying that the induced map in homology is an
isomorphism. One usually also requires that $Y$ itself is a
$l$-dimensional Poincar\'e space (using the corresponding definition
where the second space is empty) with
$\boundary[X,Y]=[Y]$, although this is not really necessary for our
applications. 
Here $C_*(\overline{X})$ is the cellular
$\Q \Gamma$-chain complex and $C^{l-*}(\overline{X},\overline{Y})$
is the dual $\Q \Gamma$-chain complex
$\hom_{\Q \Gamma}(C_{l-*}(\overline{X},\overline{Y}),\Q \Gamma)$.
Examples for Poincar\'e pairs are given by a compact connected
topological oriented manifold $X$ with boundary
$Y$ or merrily by a rational homology manifold.

The Poincar\'e duality chain map of a $4n$-dimensional
Poincar\'e pair $(X,Y)$ induces an isomorphism
$H^p(X,Y;\complexs)\to H_{4n-p}(X;\complexs)$.
If we compose the inverse with the map induced  in cohomology by the
inclusion $X\into (X,Y)$ and with the natural
isomorphism $H^p(X;\complexs)\cong H_p(X;\complexs)^*$ to the dual space 
$H_p(X;\complexs)^*$ of  $H_p(X;\complexs)$, we get in the middle 
dimension $2n$ a homomorphism
\begin{equation*}
  A: H_{2n}(X;\complexs)\to H_{2n}(X;\complexs)^*
\end{equation*}
which is selfadjoint. The signature of the (oriented)
pair $(X,Y)$ is by definition the signature of the (in general
indefinite) form $A$, i.e.~the difference of the number of positive
and negative eigenvalues of the matrix representing $A$ (after
choosing a basis for $H_{2n}(X,\complexs)$
and the dual basis for $H_{2n}(X)^*)$.

The $L^2$-signature on $(\overline{X},\overline{Y})$ is defined
similarly, but one has to replace homology by
$L^2$-homology. We get then an operator $A: H_{2n}^{(2)}(\overline{X})\to
H_{2n}^{(2)}(\overline{X})$
(using the natural isomorphism of a Hilbert space
with its dual space). The $L^2$-homology is a Hilbert module over the
von Neumann algebra $\NeumannN\Gamma$ and $A$ is a selfadjoint
bounded $\Gamma$-equivariant operator.
Hence $H_{2n}^{(2)}(\overline{X})$ splits
orthogonally into the positive part of $A$, 
the negative part of $A$ and the kernel of $A$.
The difference of the
$\NeumannN\Gamma$-dimensions of the positive part and the negative part is
by definition the \emph{$L^2$-signature}.

All this can also be reformulated in terms of cohomology instead of
homology, which is
convenient e.g.~when dealing with de Rham cohomology.

An analogue of Theorem \ref{res_approxi} for $L^2$-Betti numbers
has been proved by L\"uck \cite[Theorem 0.1]{Lueck(1994c)}.

If $X$ is a smooth closed manifold,
Atiyah's $L^2$-index theorem \cite[(1.1)]{Atiyah(1976)} shows
that the signature is multiplicative under finite coverings and that
$\sign^{(2)}(\overline{X})=\sign(X_k)/[\Gamma : \Gamma_k]$
holds for $k \ge 1$. This does not work in more general situations.
Namely, the signature is in general not multiplicative under finite coverings
neither for compact smooth manifolds with boundary
(\cite[Proposition 2.12]{Atiyah-Patodi-Singer(1975b)}
together with the Atiyah-Patodi-Singer index
theorem \cite[Theorem 4.14]{Atiyah-Patodi-Singer(1975a)})
and also not for Poincar\'e spaces $X = (X,\emptyset)$
\cite[Example 22.28]{Ranicki(1992)}, \cite[Corollary 5.4.1]{Wall(1967)}).
Our result says for these cases
that the signature is multiplicative at least approximately.
For closed topological manifolds, it is known that the signature is
multiplicative under finite coverings
\cite[Theorem 8]{Schafer(1970)}. In a companion
\cite{Lueck-Schick(2001b)} \Kommentar{Add precise reference} to this paper, we prove the following
theorem, this way apparently filling a gap in the literature:

  \begin{theorem} \label{sign^{(2)}(bar M) = sign(M)}
    Let $M$ be a closed topological manifold with normal covering
    $\overline{M} \to M$. Then
    \begin{equation*}
      \sign^{(2)}(\overline{M}) = \sign(M).
    \end{equation*}
  \end{theorem}
There, we also discuss to which extend
Theorem \ref{sign^{(2)}(bar M) = sign(M)} can be true
for Poincar\'e duality spaces $X = (X,\emptyset)$. We show
\cite{Lueck-Schick(2001b)} \Kommentar{add precise reference} that
Theorem \ref{sign^{(2)}(bar M) = sign(M)}
for Poincar\'e duality spaces $X = (X;\emptyset)$ is
implied by the $L$-theory isomorphism conjecture or by (a strong form
of) the Baum-Connes conjecture provided that $\Gamma$ is torsion-free.



Dodziuk-Mathai \cite[Theorem 0.1]{Dodziuk-Mathai(1998)} give an analog of
L\"uck's approximation theorem for $L^2$-Betti numbers to F{\o}lner
exhaustions of amenable covering spaces. Similarly, we can compute the
$L^2$-signature using a F{\o}lner exhaustion.
\begin{definition} \label{definition of exhaustion}
Let $X$ be a connected compact smooth Riemannian manifold
possibly with boundary $\boundary X$ and $\overline{X} \to X$ be
a $\Gamma$-covering for some amenable group $\Gamma$.
Let $X_1\subset X_2\subset\dots
\overline{X}$ with $\bigcup_{k\in\naturals}
X_k=\overline{X}$ be an exhaustion
of $(\overline{X},\partial \overline{X})$ by
smooth submanifolds with boundary (where
we don't make any assumptions about the intersection of $\boundary X_k$
and $\boundary \overline{X}$). Set
 $Y_k:=\boundary X_k - (\boundary X_k \cap \boundary \overline{X})$ 
(i.e.~$\boundary X_k=Y_k\cup (\boundary X_k\cap \boundary \overline{X})$). 
The exhaustion is called \emph{regular} if it
has the following properties:
\begin{enumerate}\label{regexhaust}

\item \label{amarea}
$
\area(Y_k)/\vol(X_k)\xrightarrow{k\to\infty}0;
$

\item The second fundamental forms of $\boundary X_k$ in $\overline{X}$
  and each of their covariant derivatives are
  uniformly bounded (independent of $k$);
\item The boundaries $\boundary X_k$ are uniformly collared and the
  injectivity radius of $\boundary X_k$ is uniformly bounded from
  below (always uniformly in $k$).

\end{enumerate}
\end{definition}
Regular exhaustions were introduced in
\cite[p.~152]{Dodziuk-Mathai(1997)}.
The existence of such an exhaustion is equivalent to amenability of $\Gamma$
provided that the total space $\overline{X}$ is connected.

\begin{theorem}\label{amconv}
  In the situation of Definition \ref{definition of exhaustion} we get
\begin{equation*}
  \lim_{k\to\infty} \frac{\sign(X_k,\boundary X_k)}{\vol(X_k)} =
  \frac{\sign^{(2)}(\overline{X},\overline{Y})}{\vol(X)},
\end{equation*}
where the convergence of the left hand side is part of the assertion.
\end{theorem}

The assumption that the base space $X$ is connected is necessary. 

In Theorem \ref{combconv} we give a combinatorial version of Theorem
\ref{amconv} which applies to amenable exhaustions of simplicial
homology manifolds.

For smooth manifolds with boundary, the $L^2$-signature of course is
defined in terms of the intersection pairing on $L^2$-homology. In
\cite{Lueck-Schick(2001b)} \Kommentar{Add precise reference, if
  possible} we give a proof that this coincides
with the answer predicted by the $L^2$-index theorem \cite[Theorem
1.1]{Ramachandran(1993)}. The latter paper only deals with the $L^2$-index
of certain operators, we also check the homological interpretation.

\textbf{Organization of the paper:}
 We will prove convergence of the signature for coverings in
Section \ref{sec:resconv}, and in Section \ref{sec:amenable}
the statement about amenable exhaustions.

\typeout{--------------------   Section 1 --------------------}

\section{Residual convergence of signatures}
\label{sec:resconv}

This section is devoted to the proof of Theorem \ref{res_approxi}.

Let $C_*$ be a finitely generated based free
$4n$-dimensional $\Q\Gamma$-chain complex.  Finitely generated
based free means that each chain module $C_p$ is of the shape
$\Q\Gamma^r = \oplus_{i=1}^r \Q\Gamma$ for some integer $r \ge 0$.
Its \emph{dual $\Q\Gamma$-chain complex}
$C^{4n-*}$ has as $p$-th chain module $C_{4n-*}$ and its $p$-th
differential $c^{4n-p}: C^{4n-p} \to C^{4n-(p-1)}$ is given by
$\left(c_{2d-(p-1)}: C_{2d-(p-1)} \to C_{2d-p}\right)^*$.
The adjoint $f^*: \Q\Gamma^s \to \Q\Gamma^r$ of a $\Q\Gamma$-map
$f: \Q\Gamma^r \to \Q\Gamma^s$ is
given by the matrix $A^* \in M(s,r,\Q\Gamma)$
if $f$ is given by the matrix $(A_{ij}) \in M(r,s,\Q\Gamma)$ and
$A^*_{i,j} = \overline{A_{j,i}}$ for
$\overline{\sum_{w \in \Gamma} \lambda_w \cdot w} :=
\sum_{w \in \Gamma} \lambda_w \cdot w^{-1}$. If we identify
$\hom_{\Q\Gamma}(\Q\Gamma^r,\Q\Gamma)$ with $\Q\Gamma^r$ in the obvious way,
then $f^*$ is $\hom_{\Q\Gamma}(f,\id_{\Q\Gamma})$. Given a
$\Q\Gamma$-chain map $f_*: C^{4n-*} \to C_*$, define its adjoint
$\Q\Gamma$-chain map
$f^{4n-*}: C^{4n-*} \to C_*$ in the obvious way.

Define the finitely generated $4n$-dimensional
 Hilbert $\NeumannN\Gamma$-chain complex $C_*^{(2)}$ by
 $l^2(\Gamma) \otimes_{\Q\Gamma} C_*$
and the finitely generated based free $4n$-dimensional
$\Q[\Gamma/\Gamma_k]$-chain complex $C_*[k]$ by 
$\Q[\Gamma/\Gamma_k] \otimes_{\Q\Gamma} C_*$. This applies also to
chain maps. Notice that $(C_*^{(2)})^{4n-*}$ is the same as
$(C^{4n-*})^{(2)}$ and will be denoted by $C^{4n-*}_{(2)}$
and similarly for $C_*[k]$.

Let $f_*: C^{4n-*} \to C_*$ be a $\Q\Gamma$-chain map
such that $f_*$ and its dual $f^{4n-*}$ are $\Q\Gamma$-chain homotopic.
Then both $H_{2n}^{(2)}(f_*^{(2)})$ and $H_{2n}(f_*[k])$
are selfadjoint.
Given a selfadjoint map $g: V \to V$ of Hilbert $\NeumannN\Gamma$-modules
and an interval $I \subset \reals$,
let $\chi_{I}(g)$ be the map obtained from $g$ by functional calculus
for the characteristic function $\chi_I: \reals \to \reals$ of $I$.
Define
\begin{align*}
b^{(2)}_+(g) & :=  \tr_{\NeumannN\Gamma}(\chi_{(0,\infty)}(g)); \qquad
b^{(2)}_-(g)  :=  \tr_{\NeumannN\Gamma}(\chi_{(-\infty,0)}(g));
\\
b^{(2)}(g) & :=  \dim_{\NeumannN\Gamma}(\ker(g))  =
\tr_{\NeumannN\Gamma}(\chi_{\{0\}}(g));
\\
\sign^{(2)}(g)
& :=
b^{(2)}_+(g) - b^{(2)}_-(g).
\end{align*}
If $h: W \to W$ is a selfadjoint endomorphism of a finite-dimensional
complex vector space, define analogously
\begin{align*}
b_+(h) & :=  \tr_{\complexs}(\chi_{(0,\infty)}(h));
\qquad b_-(h)  :=  \tr_{\complexs}(\chi_{(-\infty,0)}(h));
\\
b(h) & :=  \dim_{\complexs}(\ker(h)) ~ = ~
\tr_{\complexs}(\chi_{\{0\}}(h));
\\
\sign(h)
& :=
b_+(h) - b_-(h).
\end{align*}
Of course, $\sign(h)$ is the difference of the number of positive
and of negative eigenvalues of $h$ (counted with multiplicity).
Define
\begin{align*}
b_{2n\pm}^{(2)}(f_*^{(2)}) & :=
b_{\pm}^{(2)}(H_{2n}^{(2)}(f_*)); & b_{2n\pm}(f_*[k])
& :=  b_{\pm}^{(2)}(H_{2n}(f_*[k]));
\\
b_{2n}^{(2)}(f_*^{(2)}) &  :=  b^{(2)}(H_{2n}^{(2)}(f_*));&
b_{2n}(f_*[k]) & :=  b(H_{2n}(f_*[k]));
\\
\sign^{(2)}(f_*^{(2)}) & :=  \sign^{(2)}(H_{2n}^{(2)}(f_*)); &
\sign(f_*[k]) & :=  \sign(H_{2n}(f_*[k])).
\end{align*}

A classical result proved e.g.~in \cite{Lueck-Schick(2001b)}, or
(with much more information) in
\cite{Ranicki(1980aa),Ranicki(1980bb)} says
that, given a $4n$-dimensional
Poincar\'e pair
$(X,Y)$ with $\Gamma$-covering $\overline{X} \to X$, the composition of the
Poincar\'e $\Q\Gamma$-chain map
$- \cap [\overline{X},\overline{Y}]:
C^{4n-*}(\overline{X},\overline{Y}; \Q) \to
C_*(\overline{X};\Q)$ with the $\Q\Gamma$-chain map induced
by the inclusion
yields a $\Q\Gamma$-chain map
$$g_*\colon C^{4n-*}(\overline{X},\overline{Y};\Q) \to
C_*(\overline{X},\overline{Y};\Q)$$
of finitely generated
based free $4n$-dimensional $\Q\Gamma$-chain complexes such that
$g_*$ is $\Q\Gamma$-chain homotopic to $g^{4n-*}$. Define
\begin{align*}
b_{2n\pm}^{(2)}(\overline{X},\overline{Y})  &:=
b_{2n\pm}^{(2)}(g_*^{(2)}); & b_{2n\pm}(X_k,Y_k) & :=  b_{2n\pm}(g_*[k]);\\
b_{2n}^{(2)}(\overline{X},\overline{Y}) &  :=
b_{2n}^{(2)}(g_*^{(2)}); & b_{2n}(X_k,Y_k) & :=  b_{2n}(g_*[k]);
\\
\sign^{(2)}(\overline{X},\overline{Y}) & :=  \sign^{(2)}(g_*^{(2)}) ;
&
\sign(X_k,Y_k) & :=  \sign(g_*[k]).
\end{align*}

Theorem \ref{res_approxi} is an immediate consequence of
\begin{theorem} \label{the: algebraic convergence theorem}
Let $g_*\colon C^{4n-*}(\overline{X},\overline{Y};\Q)
\to C_*(\overline{X},\overline{Y};\Q)$
be the $\Q\Gamma$-chain map introduced above.
Then
\begin{align*}
b_{2n\pm}^{(2)}(g_*^{(2)}) & =  \lim_{k \to \infty}
\frac{b_{2n\pm}(g_*[k])}{[\Gamma:\Gamma_k]}.
\end{align*}
\end{theorem}

The proof of
Theorem \ref{the: algebraic convergence theorem} is split into a
sequence of lemmas.

\begin{lemma}\label{traceconv}
  Let $A:l^2(\Gamma)^n\to l^2(\Gamma)^n$ be a selfadjoint Hilbert
  $\NeumannN\Gamma$-module morphism. Let $q_j:\reals\to\reals$ be a
  sequence of measurable functions converging pointwise to the
  function $q$ such that $\abs{q_j(x)}\le C$ on the spectrum of $A$,
  where $C$ does not depend on $k$. Then
  \begin{equation*}
    \tr_{\NeumannN\Gamma}(q_j(A))\xrightarrow{j\to\infty} \tr_{\NeumannN\Gamma}(q(A)).
  \end{equation*}
\end{lemma}
\begin{proof}
  By the spectral theorem, $q_j(A)$ converges strongly to
  $q(A)$. Moreover, $\norm{q_j(A)}\le C$ for $j \in \integers$. By
  \cite[p.~34]{Dixmier(1969)} $q_j(A)$ converges ultra-strongly and
  therefore ultra-weakly to
  $q(A)$. Since $l^2(\Gamma)^n$ is a finite
  Hilbert-$\NeumannN\Gamma$-module $\1:l^2(\Gamma)^n\to l^2(\Gamma)^n$
  is of $\Gamma$-trace class. Normality of the $\Gamma$-trace implies
  the conclusion (compare \cite[Proposition 2 on p.~82]{Dixmier(1969)} or
  \cite[Theorem 2.3(4)]{Schick(2001)}).
\end{proof}

Let $\Delta_p := c_{p+1}\circ c_{p+1}^* + c_{p}^*
\circ c_p: C_p \to C_p$ be the
\emph{combinatorial Laplacian} on $\overline{X}$, where we abbreviate
$C_p:=C_p(\overline{X},\overline{A};\Q)$. Using a cellular basis of $C_p$
coming from $C_p(\overline{X},\overline{A};\integers)$
this is given by a matrix over $\integers\Gamma$. Then
$\Delta^{(2)}_p=c^{(2)}_{p+1}{c^{(2)}_{p+1}}^*+ {c_p^{(2)}}^*c_p^{(2)}: C_p^{(2)} \to C_p^{(2)}$ is the Laplacian of
$C_*^{(2)}$ and
$\Delta_p[X_k]=c_{p+1}[k]{c_{p+1}}[k]^*+c_p[k]^*c_p[k]$ is the Laplacian
on $C_p[k]$. Let $f_*: C^{4n-*}\to C_*$ be homotopic to its adjoint as
introduced in the beginning of this section.
The next lemma follows from
\cite[Lemma 2.5]{Lueck(1994c)}.

\begin{lemma}
  \label{lem: normbound}
  There is $K\ge 1$ such that for all $k \ge 1$
  \begin{equation*}
    \norm{\Delta_{2n}^{(2)}}, \norm{\Delta_{2n}[X_k]},\norm{f_{2n}^{(2)}},
    \norm{f_{2n}[k]} \le K.
  \end{equation*}
\end{lemma}

In the sequel we write
\begin{equation*}
\tr_k :=  \frac{\tr_{\Q}}{[\Gamma:\Gamma_k]};
\quad
\dim_k  :=  \frac{\dim_{\Q}}{[\Gamma:\Gamma_k]};
\quad
\sign_k  := \frac{\sign_{\Q}}{[\Gamma:\Gamma_k]},
\end{equation*}
and denote by $\pr_{2n}^{(2)}: C_{2n}^{(2)} \to C_{2n}^{(2)}$ and
$\pr_{2n}[X_k]: C_{2n}[k] \to C_{2n}[k]$ the orthogonal projection
onto the kernel of $\Delta_{2n}^{(2)}$ and $\Delta_{2n}[X_k]$.

For each $\epsilon>0$ fix a polynomial $p^{\epsilon}(x) \in \reals[x]$
with real coefficients satisfying $p^{\epsilon}(0)=1$,
$0\le p^{\epsilon}(x)\le 1+\epsilon$ for $\abs{x}\le\epsilon$ and $0\le
p(x)\le \epsilon$ for $\epsilon\le\abs{x}\le K$ (where $K$ is the
constant of Lemma \ref{lem: normbound}).

\begin{lemma}\label{balanced1}
  For each $p$ and $k$ we have
  \begin{equation*}
    \dim_k C_p(X_k) = \dim_{\NeumannN\Gamma} C_p^{(2)}(\overline X),
  \end{equation*}
  and hence in particular
  \begin{equation*}
   \lim_{k\to\infty} \dim_k C_p(X_k) = \dim_{\NeumannN\Gamma} C_p^{(2)}(\overline X).
  \end{equation*}
\end{lemma}
\begin{proof}
  For every $k$, $\dim_k C_p(X_k)$ is equal to the number of $p$-cells
  in $X$, and the same is true for $\dim_{\NeumannN\Gamma} C_p^{(2)}(\overline X)$.
\end{proof}

\begin{lemma}\label{convlemma}
For $\Q\Gamma$-linear maps $h_1,\dots,h_d: \Q\Gamma^r \to \Q\Gamma^r$
and a polynomial $p(x_1,\dots,x_d)$ in non-commuting variables
$x_1,\dots,x_d$ we have
\begin{equation*}
    \tr_{\NeumannN\Gamma}(p(h_1^{(2)},\dots, h_d^{(2)})) =
    \lim_{k\to\infty}tr_k\left(p(h_1[k],\dots, h_d[k])\right) .
\end{equation*}
\end{lemma}
\begin{proof}
  By linearity it suffices to prove this for monomials $p=x_{i_1}\dots
  x_{i_d}$, and since the $h_j$ are not assumed to be
  different, without loss of generality we can assume $p=x_1\dots x_d$.
  The proof of \cite[Lemma 2.6]{Lueck(1994c)} applies and shows that
  there is $L>0$ such that $\tr_{\NeumannN\Gamma}(h_1^{(2)}\circ\dots\circ
  h_d^{(2)})=\tr_k(h_1[k]\circ\dots \circ h_d[k])$ for $k \ge L$.
\end{proof}
The lemma is formulated in a way that it can be applied if the
assignment $h\to h[k]$ is not a homomorphism. This is unnecessary
here, but will be needed in Section \ref{sec:amenable}.

\begin{lemma}\label{lem:smallEV1}
  There is a
constant $C_1>0$ (independent of $k$) such that for
$0 < \epsilon < 1$ and $k \ge 1$
\begin{equation}
\tr_k\left(\chi_{(0,\epsilon]}(\Delta_{2n}[X_k])\right) \le
\frac{C_1}{-\ln(\epsilon)}.
\end{equation}
\end{lemma}
\begin{proof}
   This is part of \cite[Lemma 2.8]{Lueck(1994c)}.
\end{proof}

\begin{lemma}\label{lem: speccontrol}
  We find a constant $C>0$ (independent of $k$)
  such that for all $k \ge 1$ and $0 <\epsilon < 1$
  \begin{equation*}
  0 \le     \tr_k \left(
  \abs{p^\epsilon(\Delta_{2n}[X_k])-\pr_{2n}[X_k]}\right) \le
   C\cdot \epsilon +  \frac{C}{-\ln(\epsilon)}.
  \end{equation*}
  Moreover
  \begin{equation*}
    \lim_{\epsilon \to 0}
    \tr_{\NeumannN\Gamma}
    \left(\left|p^{\epsilon}(\Delta_{2n}^{(2)})-\pr_{2n}^{(2)}\right|\right) = 0.
  \end{equation*}
\end{lemma}

\begin{proof}
First observe that by our construction
$p^{\epsilon}(\Delta_{2n}[X_k])-\pr_{2n}[X_k]$
is non-negative since $0\le p^{\epsilon}-\chi_{\{0\}}$ on
 the spectrum. We also have
 $p^\epsilon-\chi_{\{0\}} \le \epsilon + \chi_{(0,\epsilon]}$
 on the spectrum of the operators. Since the trace is positive, we get
  \begin{eqnarray*}
   0\le  \tr_k(p^{\epsilon}(\Delta_{2n}[X_k])-\pr_{2n}[X_k])  & \le &
    \epsilon  \tr_k(\id_{C_{2n}[k]}) +
    \tr_k(\chi_{(0,\epsilon]}(\Delta_{2n}[X_k])).
  \end{eqnarray*}
Now the first inequality follows from Lemma \ref{balanced1} and
Lemma \ref{lem:smallEV1}. The second one follows from
  \begin{eqnarray*}
    \tr_{\NeumannN\Gamma} \left(p^{\epsilon}
    (\Delta_{2n}^{(2)})-\pr_{2n}^{(2)}\right)
    & \le &
    \epsilon  \tr_{\NeumannN\Gamma}\left(\id_{C_{2n}^{(2)}}\right) +
    \tr_{\NeumannN\Gamma}\left(\chi_{(0,\epsilon]}(\Delta_{2n}^{(2)})\right)
  \end{eqnarray*}
  and the fact
  that because of Lemma \ref{traceconv}
  $\lim_{\epsilon \to 0}
  \tr_{\NeumannN\Gamma}\left(\chi_{(0,\epsilon]}(\Delta_{2n}^{(2)})\right)
  =0$.
\end{proof}

We also cite the following result \cite[Theorem 2.3]{Lueck(1994c)}:
\begin{theorem}
  \label{betticonv1}
  The normalized sequence of Betti numbers converges, i.e.~for each $p$
  \begin{equation*}
   \lim_{k\to\infty}    \dim_k(\ker(\Delta_p[X_k])) = \dim_{\NeumannN\Gamma}
   \ker(\Delta^{(2)}_p).
  \end{equation*}
\end{theorem}
We want to approximate $\chi_{(a,b)}$ by polynomials. Next we check
that for a fixed polynomial we can replace $\pr_{2n}[X_k]$ in the argument by
$p^\epsilon(\Delta_{2n}[X_k])$.

\begin{lemma}\label{lem: uniformconv}
  Fix a polynomial $q\in\reals [x]$. Then
  we find a constant $D>0$ (independent of $k$)
  such that for all $k \ge 1$ and $0 < \epsilon < 1$
  \begin{multline*}
    |
    \tr_k\left(q\left(p^\epsilon(\Delta_{2n}[X_k]) \circ f_{2n}[k]
    \circ p^\epsilon(\Delta_{2n}[X_k])\right)\right) -\\
    \tr_k\left(q\left(\pr_{2n}[X_k] \circ f_{2n}[k]
    \circ \pr_{2n}[X_k]\right)\right)
    | 
     \le  D\cdot \epsilon + \frac{D}{-\ln(\epsilon)}.
  \end{multline*}
  Moreover, we have
  \begin{equation*}
    \lim_{\epsilon\to 0}
    \tr_{\NeumannN\Gamma}\left(q\left(p^\epsilon(\Delta_{2n}^{(2)}) \circ
f_{2n}^{(2)} \circ p^\epsilon(\Delta_{2n}^{(2)})
    \right)\right)  =
    \tr_{\NeumannN\Gamma}\left(q\left(\pr_{2n}^{(2)}
    \circ f_{2n}^{(2)} \circ \pr_{2n}^{(2)}\right)\right).
  \end{equation*}
\end{lemma}
\begin{proof}
  By linearity it suffices to prove the statement
  for all monomials $q(x)=x^n$. Obviously it suffices to consider
  $n \ge 1$. In the sequel we abbreviate
  $x = p^\epsilon(\Delta_{2n}[X_k])$, $f = f_{2n}[k]$ and
   $y=  \pr_{2n}[X_k]$.   Notice that $\norm{x} \le (1+ \epsilon)$,
  $\norm{f} \le K$ and $\norm{y}\le 1$ holds for
  the constant $K$ appearing in Lemma \ref{lem: normbound}.
  We estimate using the
  trace property $\tr(AB)=\tr(BA)$ and the trace
  estimate $\abs{\tr(AB)}\le  \norm{A} \cdot \tr(\abs{B})$ (which
  also holds for the normalized traces $\tr_k$ and for $\tr_{\NeumannN\Gamma}$
  by \cite[p.~106]{Dixmier(1969)} since all the traces we are
  considering are normal),
  \begin{eqnarray*}
  \lefteqn{\left|\tr_k\left((p^\epsilon(\Delta_{2n}[X_k]) \circ f_{2n}[k]
  \circ p^\epsilon(\Delta_{2n}[X_k]))^n\right)\right.} & &\\
 & & -
   \left. \tr_k\left((\pr_{2n}[X_k] \circ f_{2n}[k]
    \circ \pr_{2n}[X_k])^n\right)\right|
  \\ & = &
  \left|\tr_k\left(xfxxfx\ldots xfx - yfyyfy\ldots yfy\right)\right|
    \\ & = &
  \left|\tr_k\left((x-y)fxxfx \ldots xfx + yf(x-y)xfx \ldots xfx
  \right.\right.
   \\ & = &
  \hspace{20mm} +
  \left. \left.yfy(x-y)fxxfx\ldots fx +
  \ldots + yfyyfy \ldots yf(x-y)\right)\right|
   \\ & \le &
   2n \cdot (1+\epsilon)^{2n-1} \cdot K^n \cdot \tr(\abs{x-y})
   \\ & = &
    2n \cdot (1+\epsilon)^{2n-1} \cdot K^n \cdot
   \tr_k\left(\left|p^\epsilon(\Delta_{2n}[X_k]) -  \pr_{2n}[X_k]\right|\right).
\end{eqnarray*}
A similar estimate holds in the $L^2$-case. The claim follows from
Lemma \ref{lem: speccontrol}.
\end{proof}

\begin{lemma} \label{lem: lim for q}
Fix a polynomial $q(x) \in \reals[x]$. Then
$$\lim_{k \to \infty} \tr_k\left(q\left(\pr_{2n}[X_k]
\circ f_{2n}[k]\circ \pr_{2n}[X_k]\right)\right)
~ = ~ \tr_{\NeumannN\Gamma}\left(q\left(\pr_{2n}^{(2)} \circ
f_{2n}^{(2)} \circ \pr_{2n}^{(2)}\right)\right).$$
\end{lemma}
\begin{proof}
Fix $\delta >0$. By Lemma \ref{lem: uniformconv} we find
$\epsilon > 0$   such that for all $k \ge 1$
\begin{multline*}
\bigl|\tr_k\left(q\left(p^\epsilon(\Delta_{2n}[X_k]) \circ f_{2n}[k]
    \circ p^\epsilon(\Delta_{2n}[X_k])\right)\right) -\\
    \tr_k\left(q\left(\pr_{2n}[X_k] \circ f_{2n}[k]
    \circ \pr_{2n}[X_k] \right)\right)\bigr|
\le  \delta/3;
\end{multline*}
\begin{multline*}
    \bigl|\tr_{\NeumannN\Gamma}\left(q\left(p^\epsilon(\Delta_{2n}^{(2)})
    \circ f_{2n}^{(2)} \circ p^\epsilon(\Delta_{2n}^{(2)})\right)\right)  -\\
    \tr_{\NeumannN\Gamma}\left(q\left(\pr_{2n}^{(2)}
    \circ f_{2n}^{(2)} \circ \pr_{2n}^{(2)}
    \right)\right)\bigr|
    \le  \delta/3.
\end{multline*}
Hence it suffices to show for each fixed $\epsilon$
\begin{multline*}
 \lim_{k \to \infty}  \tr_k
 \bigl( q \bigl(p^\epsilon(c_{p+1}[k]c_{p+1}[k]^*+\\
   c_p[k]^* c_p[k] ) \circ f_{2n}[k]
\circ p^\epsilon(c_{p+1}[k]c_{p+1}[k]^*+c_p[k]^* c_p[k])\bigr)\bigr)\\
     =  
    \tr_{\NeumannN\Gamma}\bigl(q
    \left(p^\epsilon(c_{p+1}^{(2)}{c_{p+1}^{(2)}}^*+{c_p^{(2)}}^* c_p^{(2)}) \circ f_{2n}^{(2)}
     \circ p^\epsilon(c_{p+1}^{(2)}{c_{p+1}^{(2)}}^* +{c_p^{(2)}}^*
     c_p^{(2)})\right)\bigr). 
\end{multline*}
Since $q$ and $p^{\epsilon}$ are fixed, we deal with a fixed polynomial
expression in $c_p$, $c_p^*$, $c_{p+1}$, $c_{p+1}^*$, and
$f_{2n}$. Therefore the last claim
follows from Lemma \ref{convlemma}.
This finishes the proof of Lemma \ref{lem: lim for q}.
\end{proof}

\begin{lemma} \label{lem: liminf-estimate}
We have for $a,b \in\reals $ with  $a<b$
$$\tr_{\NeumannN\Gamma}\left(\chi_{(a,b)}
\left(H_p^{(2)}(f_*^{(2)})\right)\right)
~ \le ~ \liminf_{k\to\infty}
\tr_k\left(\chi_{(a,b)}\left(H_p(f_*[k])\right)\right).$$
\end{lemma}
\begin{proof}
We approximate $\chi_{(a,b)}$ by polynomials. Namely,
for $0 < \epsilon < (b-a)/2$ and $K$ as above
let $q^\epsilon\in\reals[x]$ be a polynomial with
\begin{eqnarray*}
-1 \le  q^\epsilon(x)\le \chi_{(a,b)}(x)
& &  \mbox{ for }
     \abs{x}\le K;\\
q^\epsilon(x) \ge \chi_{(a,b)}(x) - \epsilon
    & & \mbox{ for } x\in [-K, a]\cup [ a+\epsilon, b-\epsilon]\cup [ b, K].
\end{eqnarray*}
Under the identification of $\im(\pr_{2n}[X_k])$ and $H_p(C_*[k])$
coming from the (combinatorial) Hodge-de Rham theorem
the operator
$\pr_{2n}[X_k]\circ  f_{2n}[k] \circ \pr_{2n}[X_k]$  restricted to
$\im(\pr_{2n}[X_k])$ becomes $H_p(f_*[k])$ which is selfadjoint
because of $f_* \simeq f^{4n-*}$. Hence
$\pr_{2n}[X_k]\circ  f_{2n}[k] \circ \pr_{2n}[X_k]$ and also the operator
$q^{\epsilon}(\pr_{2n}[X_k]\circ  f_{2n}[k] \circ \pr_{2n}[X_k])$
are selfadjoint.
The same is true on
the $L^2$-level and we conclude
\begin{eqnarray}
\hspace{-6mm} \tr_k\left(\chi_{(a,b)}\left(\pr_{2n}[X_k]\circ
f_{2n}[k]\circ \pr_{2n}[X_k]\right)\right)
& = &
\tr_k\left(\chi_{(a,b)}\left(H_p(f_*[k])\right)\right); \label{eqn 1}
\\
\hspace{-6mm} \tr_{\NeumannN\Gamma}\left(\chi_{(a,b)}\left(\pr_{2n}^{(2)}\circ
f_{2n}^{(2)}\circ \pr_{2n}^{(2)}\right)\right)
& = &
\tr_{\NeumannN\Gamma}
\left(\chi_{(a,b)}\left(H_p^{(2)}(f_*^{(2)})\right)\right).
\label{eqn 2}
\end{eqnarray}
Positivity of the trace  and $q^{\epsilon}(x) \le \chi_{(a,b)}(x)$
for all $x$ in the spectrum of
$\pr_{2n}[X_k] \circ f_{2n}[k]\circ \pr_{2n}[X_k]$ implies
\begin{equation*}
  \tr_k\left(q^{\epsilon}\left(\pr_{2n}[X_k]
\circ f_{2n}[k] \circ \pr_{2n}[X_k] \right)\right) \le
  \tr_k\left(\chi_{(a,b)}\left(\pr_{2n}[X_k]
\circ  f_{2n}[k]\circ \pr_{2n}[X_k]\right)\right).
\end{equation*}
Note that for fixed $q^\epsilon$ the left hand side
converges for $k \to \infty$ by Lemma \ref{lem: lim for q}.
For the right hand side this is not
clear, but in any case we get
\begin{multline} \label{eqn 3}
  \tr_{\NeumannN\Gamma}\left(q^{\epsilon}
  \left(\pr_{2n}^{(2)} \circ f_{2n}^{(2)}
    \circ \pr_{2n}^{(2)}\right)\right)\\
  \le
  \liminf_{n\to\infty}
  \tr_k\left(\chi_{(a,b)}\left(\pr_{2n}[X_k]
  \circ  f_{2n}[k]\circ \pr_{2n}[X_k]\right)\right).
\end{multline}

On the spectrum of the operator in question, the functions
$q^\epsilon$ are uniformly bounded and converge
pointwise to
$\chi_{(a,b)}$ if $\epsilon\to 0$. By Lemma \ref{traceconv}
\begin{equation*}
  \lim_{\epsilon\to 0}
  \tr_{\NeumannN\Gamma}\left(q^{\epsilon}
  \left(\pr_{2n}^{(2)} \circ f_{2n}^{(2)}
  \circ \pr_{2n}^{(2)}\right) \right)
   =
   \tr_{\NeumannN\Gamma}\left(\chi_{(a,b)}\left(\pr_{2n}^{(2)}
\circ f_{2n}^{(2)} \circ \pr_{2n}^{(2)} \right)\right).
\end{equation*}
Since inequality \eqref{eqn 3} holds for arbitrary $\epsilon>0$, we
conclude
\begin{multline*}
\tr_{\NeumannN\Gamma}\left(\chi_{(a,b)}\left(\pr_{2n}^{(2)}
\circ f_{2n}^{(2)}\circ \pr_{2n}^{(2)} \right)\right)\\
\le \liminf_{k\to\infty}
  \tr_k\left(\chi_{(a,b)}\left(\pr_{2n}[X_k]\circ
 f_{2n}[k]\circ \pr_{2n}[X_k]\right)\right).
\end{multline*}
Now the claim follows from
\eqref{eqn 1} and \eqref{eqn 2}.
\end{proof}

\begin{lemma} \label{lem: conv of homology kernels}
Let $f_*: C_* \to D_*$ be a $\rationals\Gamma$-chain
map of finitely generated
based free $\rationals\Gamma$-chain complexes. Then we get for all $p$
$$\lim_{k \to \infty} \dim_k\left(\ker\left(H_p(f_*[k]\right))\right)
~ = ~
\dim_{\NeumannN\Gamma}\left(\ker\left(H_p^{(2)}(f_*^{(2)})\right)\right).
$$
\end{lemma}
\begin{proof} We can assume without loss of
generality that $C_*$ and $D_*$ are $(p+1)$-dimensional.
Consider the long exact sequence of $\Q \Gamma$-chain complexes
$0 \to D_* \to \cone(f_*)_* \to \Sigma C_* \to 0$,
where $\cone(f_*)_*$ is the mapping cone of
$f_*$ and $\Sigma C_*$ the suspension of $C_*$. It is a split exact sequence in each
dimension and thus remains exact after applying
 $l^2(\Gamma) \otimes_{\Q \Gamma} -$.
The weakly exact long homology sequence yields
a weakly exact sequence of Hilbert $\NeumannN(\Gamma)$-modules
$$0 \to
H_{p+2}^{(2)}(\cone(f_*)_*^{(2)})
\to H_{p+1}^{(2)}(C_*^{(2)})
\xrightarrow{H_{p+1}^{(2)}(f_*^{(2)})}  H_{p+1}^{(2)}(D_*^{(2)})$$
$$ \to H_{p+1}^{(2)}(\cone(f_*)^{(2)}_*))
\to \ker(H_p(f_*^{(2)})) \to 0.$$
This implies
\begin{eqnarray}
\lefteqn{\dim_{\NeumannN\Gamma}\left(\ker(H_p(f_*^{(2)}))\right)}
&  & \nonumber
\\
& = & \dim_{\NeumannN\Gamma}\left(H_{p+1}^{(2)}(\cone(f_*)^{(2)}_*)\right) -
\dim_{\NeumannN\Gamma}\left(H_{p+1}^{(2)}(D_*^{(2)})\right) \nonumber
\\
& & +
\dim_{\NeumannN\Gamma}\left(H_{p+1}^{(2)}(C_*^{(2)}) \right) -
\dim_{\NeumannN\Gamma}\left(H_{p+2}^{(2)}(\cone(f_*)_*^{(2)})\right).
\label{eqn 11}
\end{eqnarray}
Analogously we get
\begin{eqnarray}
\lefteqn{\dim_{k}\left(\ker(H_p(f_*[k]))\right)} &  & \nonumber
\\
& = & \dim_{k}\left(H_{p+1}(\cone(f_*[k])_*)\right) -
\dim_{k}\left(H_{p+1}(D_*[k])\right) \nonumber
\\
& & +
\dim_{k}\left(H_{p+1}(C_*[k]) \right) -
\dim_{k}\left(H_{p+2}(\cone(f_*[k])_*)\right). \label{eqn 12}
\end{eqnarray}
We conclude from Theorem \ref{betticonv1}
\begin{eqnarray}
\dim_{\NeumannN\Gamma}\left(H_{p+1}^{(2)}(\cone(f_*)^{(2)}_*)\right)
& = &
\lim_{k \to \infty} \dim_{k}\left(H_{p+1}(\cone(f_*[k])_*)\right);
\label{eqn 13}
\\
\dim_{\NeumannN\Gamma}\left(H_{p+1}^{(2)}(D_*^{(2)})\right)
& = &
\lim_{k \to \infty} \dim_{k}\left(H_{p+1}(D_*[k])\right);
\label{eqn 14}
\\
\dim_{\NeumannN\Gamma}\left(H_{p+1}^{(2)}(C_*^{(2)}) \right)
& = &
\lim_{k \to \infty} \dim_{k}\left(H_{p+1}(C_*[k]) \right);
\label{eqn 15}
\\
\dim_{\NeumannN\Gamma}\left(H_{p+2}^{(2)}(\cone(f_*)_*^{(2)})\right)
& = &
\lim_{k \to \infty} \dim_{k}\left(H_{p+2}(\cone(f_*[k])_*)\right).
\label{eqn 16}
\end{eqnarray}

Now the claim follows from equations  \eqref{eqn 11}--\eqref{eqn 16}.
\end{proof}

Now we are ready to prove Theorem
\ref{the: algebraic convergence theorem}.
\begin{proof}
We get from Lemma \ref{lem: liminf-estimate} and
 Lemma \ref{lem: conv of homology
  kernels}
\begin{align*}
b_{2n+}^{(2)}(g_*^{(2)}) & \le  \liminf_{k \to \infty}
\frac{b_{2n+}(g_*[k])}{[\Gamma:\Gamma_k]};&
b_{2n-}^{(2)}(g_*^{(2)}) & \le  \liminf_{k \to \infty}
\frac{b_{2n-}(g_*[k])}{[\Gamma:\Gamma_k]};\\
b_{p}^{(2)}(g_*^{(2)}) & =  \lim_{k \to \infty}
\frac{b_{p}(g_*[k])}{[\Gamma:\Gamma_k]}.
\end{align*}
Since
\begin{eqnarray}
b_{2n+}^{(2)}(g_*^{(2)}) + b_{2n-}^{(2)}(g_*^{(2)}) +
b_{2n}^{(2)}(g_*^{(2)}) & = &
\dim_{\NeumannN\Gamma}\left(C_{2n}^{(2)}\right);\nonumber
\\
\frac{b_{2n+}(g_*[k])}{\abs{\Gamma/\Gamma_k}} +
\frac{b_{2n-}(g_*[k])}{\abs{\Gamma/\Gamma_k}} + 
\frac{b_{2n}(g_*[k])}{\abs{\Gamma/\Gamma_k}} & = & \dim_{k}(C_{2n}[k]);\nonumber
\\
\dim_{\NeumannN\Gamma}\left(C_{2n}^{(2)}\right) = \lim_{k\to\infty} \dim_k(C_{2n}[k])
& = &
\frac{\dim_{\Q}(C_{2n}[k])}{[\Gamma:\Gamma_k]}, \nonumber
\end{eqnarray}
Theorem \ref{the: algebraic convergence theorem} and thus 
Theorem \ref{res_approxi} follow from Lemma \ref{balanced1}.
\end{proof}

\begin{remark}\label{etaremark}
  Theorem \ref{res_approxi} can be applied to a
  $4n$-dimensional Riemannian manifold $X$ with boundary $Y$. In this case, the
  Atiyah-Patodi-Singer  theorem \cite[Theorem
  4.14]{Atiyah-Patodi-Singer(1975a)} and
  \cite[(0.9)]{Cheeger-Gromov(1985a)} and the $L^2$-signature theorem
  of \cite{Lueck-Schick(2001b)} \Kommentar{Add more precise reference}
imply
\begin{equation*}
  \begin{split}
    \frac{\sign(X_k,\boundary X_k)}{\vol(X_k)} =& \frac{1}{\vol(X_k)}\cdot \int_{X_k}
    L(X_k) +
    \frac{\eta(\boundary X_k)}{\vol(X_k)} +
    \frac{1}{\vol(X_k)} \cdot \int_{\boundary X_k} \Pi_L(\boundary X_k),\\
 \frac{\sign^{(2)}(\overline{X},\overline{\boundary X})}{\vol(X)} = &
    \frac{1}{\vol(X)}\cdot\int_{X} L(X)
  +  \frac{\eta^{(2)}(\boundary\overline{X})}{\vol(X)} +
  \frac{1}{\vol(X)}\cdot  \int_{\boundary X} \Pi_L(\boundary X).
\end{split}
\end{equation*}
Here $L(X_k)$ and $L(X)$ denote the Hirzebruch $L$-polynomial, and
$\Pi_L(\boundary X_k)$ and $\Pi_L(\boundary X)$ are a local correction terms
which arises because the metric
is not a product near the boundary. Being local expressions, the first 
and the third summand does not depend on $k$. It follows that the sequence of
$\eta$-invariants converges. In fact, even without the assumption that
$Y^{4n-1}$ is a boundary of a suitable manifold $X$, in
\cite[Theorem 3.12]{Vaillant(1997)} it is proved
    \begin{equation*}
      \lim_{k\to\infty} \frac{\eta(Y_k)}{[\Gamma : \Gamma_k]} =
      \eta^{(2)}(\overline{Y}) .
    \end{equation*}
Key ingredients are on the one hand the analysis of Cheeger-Gromov in
\cite[Section
7]{Cheeger-Gromov(1985a)} of the formulas \eqref{def of eta} and
\eqref{def of eta^{(2)}} (which holds for operators different from
the signature operator). We present similar considerations in Section
\ref{sec:analytic}.  The second key ingredient is L\"uck's
approximation result for
$L^2$-Betti numbers \cite[Theorem 0.1]{Lueck(1994c)} (which is special
to the Laplacian, the
square of the signature operator).
\end{remark}

\begin{remark}
  The normalized signatures
  $\frac{\sign(X_k,Y_k)}{\abs{\Gamma/\Gamma_k}}$ are the
  $L^2$-signatures $\sign^{(2)}(X_k,Y_k)$ of the $\Gamma/\Gamma_k$-coverings
  $(X_k,Y_k)\to (X,Y)$. With
  this reformulation, one may ask whether Theorem \ref{res_approxi}
  holds if $\Gamma/\Gamma_k$ is not necessarily finite. 
  
This is indeed the case if the groups $\Gamma/\Gamma_k$ belong to a
large class of groups $\RAgroups$   defined in \cite[Definition
1.11]{Schick(2001b)}. 

  The
  corresponding question for $L^2$-Betti numbers is answered
  affirmatively in \cite[Theorem 6.9]{Schick(2001b)} whenever
  $\Gamma/\Gamma_k \in\RAgroups$. 
  As just mentioned, Theorem
  \ref{res_approxi} extends to this situation as well, and the proof
  we have given is formally unchanged, using the generalization of
  Lemma \ref{lem: normbound} and Lemma \ref{convlemma} given in
  \cite[Lemma 5.5 and 5.6]{Schick(2001b)}. It only remains to
  establish Lemma \ref{lem:smallEV1}, which is not done in
  \cite{Schick(2001b)}. We do this in the following Lemma
  \ref{detcontrol}, which 
  applies because of \cite[6.9]{Schick(2001b)} and because of Lemma
  \ref{lem: normbound}. 
\end{remark}

\begin{lemma}\label{detcontrol} If $\norm{\Delta[X_k]}\le K$ and
    \begin{equation}
    \label{eq:detest}
    \ln\Det'_{(2)}(\Delta[X_k]):= \int_{0^+}^\infty \ln(\lambda)
    \;dF_{\Delta[X_k]}(\lambda) \ge 0
  \end{equation}
then 
\begin{equation}\label{eq:ast}
    \tr_k(\chi_{(0,\epsilon]}(\Delta[X_k])) \le \frac{d\cdot
    \ln(K)}{-\ln(\epsilon)}.
  \end{equation}
 Here $F_{\Delta[X_k]}(\lambda)$
  is the spectral density function of the operator $\Delta[X_k]$
  computed using $\dim_k$ instead of $\dim_\complexs$, and
  $d=F_{\Delta[k]}(K)$ is the number of rows (and columns) of the
  matrix $\Delta$.
\end{lemma}
\begin{proof}
  We argue as follows (with $F:=F_{\Delta[X_k]}$):
  \begin{equation*}
  \begin{split}
    \int_{0^+}^\infty \ln(\lambda)\; dF(\lambda) & = \int_{0^+}^\epsilon
    \ln(\lambda)\; dF(\lambda) +\int_\epsilon^{\norm{\Delta[X_k]}}
      \ln(\lambda)\;dF(\lambda) \\
    & \le \ln(\epsilon)\underbrace
    {( F(\epsilon)-F(0))}_{=\tr_k(\chi_{(0,\epsilon]}(\Delta[X_k]))} +
    \ln(\norm{\Delta[X_k]}) F(\norm{\Delta[X_k]}).
  \end{split}
\end{equation*}
For $0<\epsilon<1$, using the bound $\norm{\Delta[X_k]}\le K$ of the
  generalization of Lemma \ref{lem: normbound}, Inequality
  \eqref{eq:detest} immediately gives \eqref{eq:ast}.
\end{proof}

\typeout{--------------------   Section 3 --------------------}

\section{Amenable convergence of signatures}
\label{sec:amenable}

\subsection{Analytic version}\label{sec:analytic}

In this subsection we want to prove
Theorem \ref{amconv}. We will use the following notion of manifold with bounded geometry
(compare e.g.~\cite[Definition 2.24]{Lueck-Schick(1998)}).
\begin{definition}
  A Riemannian manifold $(M,g)$ (the boundary may or may not be empty)
  is called a manifold
of {\em bounded geometry} if {\em bounded geometry constants}
  $C_q$ for $q\in M$ and
 $R_I, R_C>0$ exist, so that the following holds:
 \begin{enumerate}
\item\label{collar} The geodesic flow of the unit inward normal field induces a
diffeomorphism of $[0,2R_C)\times \boundary M$ onto its image,
 the geodesic collar;
\item
 For $x \in M$ with $d(x,\boundary M)>R_C/2$
 the exponential map $T_xM\to M$ is
a diffeomorphism on $B_{R_I}(0)$;
\item The injectivity radius of $\boundary M$ is bigger than $R_I$;
\item For every $q\in M$ we have
  $\abs{\nabla^i R}\le C_k$ and
  $\abs{\nabla_{\partial}^i l}\le C_l$ for $0\le i\le q$, where $R$ is the
  curvature tensor of $M$, $l$ the second fundamental form tensor of
  $\boundary M$, and $\nabla^i$ and $\nabla_{\partial}^i$ are the
  covariant derivatives of  $M$ and $\boundary M$.
 \end{enumerate}
\end{definition}
By \cite[Theorem 2.4]{Schick(2001a)} this
is equivalent to  \cite[Definition 2.24]{Lueck-Schick(1998)}.

Every compact manifold, or more generally every covering of a compact manifold, is a manifold with bounded geometry.

We now repeat a few well known facts about manifolds of bounded geometry.
\begin{proposition}\label{bettibound}
  Let $M$ be a compact smooth Riemannian manifold.
  There is a  constant $A>0$, depending only on the bounded geometry
  constants and the dimension of $M$, such that
  \begin{eqnarray*}
     \abs{\,\exp(-t\Delta_p(M))(x,x)} & \le & A \hspace{20mm}
     \mbox{ for } t\ge 1, x \in M;
     \\
     b_p(M)  & \le & A\vol(M);
    \\
    b_p(M,\boundary M) & \le & A\vol(M),
   \end{eqnarray*}
where the Laplacian can be taken with either relative
or absolute boundary conditions.
\end{proposition}
\begin{proof}
 The first inequality is proved in
  \cite[Theorem 2.35]{Lueck-Schick(1998)}. The claim for the Betti numbers
   is a consequence of the fact that the Betti number $ b_p(M)$
   or $b_p(M,\boundary M)$ can be written as
   $ \lim_{t\to\infty} \int_M \tr_x\exp(-t\Delta_p(M))(x,x)\;dx$
   for the Laplacian with absolute or relative boundary
  conditions, respectively. 
 \end{proof}

\begin{theorem}\label{localize}
  Let $M,N$ be Riemannian manifolds without boundary which are of bounded
  geometry and with a fixed set of bounded geometry constants. Let $U$
  be an open subset of $M$ which is isometric to a subset of $N$ (which
  we identify with $U$). For $R>0$ set 
  \begin{equation*}
U_R:=\{x\in U|\; d(x,M-U)\ge
  R\text{ and }d(x,N-U)\ge R\}.
\end{equation*}
Let $D[M]$ and $D[N]$ be the
  (tangential) signature operators 
  on $M$ and $N$, respectively; and similarly
  $\Delta[M]$ and $\Delta[N]$ the Laplacian (on differential forms). 
  Let $e^{-t\Delta}(x,y)$ and $De^{-tD^2}(x,y)$ be the integral
  kernels (which are smooth) of the operators $e^{-t\Delta}$ and
  $De^{-tD^2}$. Then there are constants $C_1,C_2>0$ which depend only
  on the dimension and the given bounded geometry constants such that
  for $t > 0$, $x \in U_R$ and $p \ge 0$
  \begin{align}\label{Delta}
    \abs{e^{-t\Delta[M]}(x,x)-e^{-t\Delta[N]}(x,x)} & \le C_1 \cdot
    e^{-R^2C_2/t};
    \\\label{sgn}
    \abs{D[M]e^{-tD[M]^2}(x,x) - D[N]e^{-tD[N]^2}(x,x)} & \le C_1\cdot
    e^{-R^2C_2/t}.
  \end{align}
\end{theorem}
\begin{proof}
  This follows by a standard argument of Cheeger-Gromov-Taylor
  \cite{Cheeger-Gromov-Taylor(1982)} from unit propagation speed and
  local elliptic estimates (here the bounded geometry constants come
  in). A detailed account is given in the proof of
  \cite[Theorem 2.26]{Lueck-Schick(1998)} which yields immediately
  \eqref{Delta}. Replacing $\sqrt{\Delta}$ by $D$ (which is possible
  since we are looking for manifolds without boundary, so that we do
  not have to worry about the non-locality of boundary conditions and
  therefore have unit propagation speed for $D$, too), the proof also
  applies to the tangential signature operator to give \eqref{sgn}.
\end{proof}

\begin{proposition}\label{smallt}
  Let $M^m$ be a manifold of bounded geometry with fixed
  bounded geometry constants and with $\boundary M=\emptyset$. Let $D$
  be the (tangential) signature operator on $M$. Then there is a function
  $A: [0,\infty) \to (0,\infty)$
  which depends only on the bounded geometry
  constants and the dimension $m$, such that for $T \ge 0$
  \begin{equation*}
    \abs{\tr_x\left(D e^{-tD^2}(x,x)\right)}  \le A(T)\cdot t^{1/2}
\hspace{10mm} \mbox{ for }
    0\le t\le T, x\in M.
  \end{equation*}
\end{proposition}
\begin{proof}

One can use the proof of
  \cite[Lemma 3.1.1 on p.~ 324]{Ramachandran(1993)} (where a slightly different
  statement is proved). The proposition is also implicit in
  \cite[Proof of Theorem 0.1 on p.~140]{Cheeger-Gromov(1985a)}. The proof uses
  the cancellation of the coefficients of
  negative powers of $t$ in the local asymptotic expansion due to
  Bismut and Freed \cite[Theorem 2.4]{Bismut-Freed(1986b)}
  and a localization
  argument based on elliptic estimates
  (here the local geometry comes in), together with the finite propagation
  speed method of Cheeger-Gromov-Taylor
  \cite{Cheeger-Gromov-Taylor(1982)}.
\end{proof}

We fix the following notation.

\begin{notation} \label{not: n_k and N_k}
In the situation of Definition \ref{definition of exhaustion} put for $r \ge 0$
$$U_r(Y_k):=\{x\in \overline{X};\;d(x,Y_k)\le r\},$$
where two points $y,z \in \overline{X}$ 
have distance $d(y,z) = d$ if there is a geodesic
of length $d$ in $\overline{X}$ joining $y$ and $z$ and $d = \infty$ 
if there is no such geodesic.
In particular $d(y,z) < \infty$ implies 
that $y$ and $z$ lie in the same path component of $\overline{X}$.
Let $\mathcal{F}$ be a (compact) connected simplicial fundamental
domain for $X$ in $\overline{X}$ such 
that $\mathcal{F} \cap \boundary \overline{X}$
is a fundamental domain for $\boundary X$. (We can construct $\mathcal{F}$ 
as a union of lifts of the top-dimensional 
simplices in a smooth triangulation of $X$
and achieve $\mathcal{F}$ to be connected, 
since $X$ is connected by assumption.)
For $r\ge 0$ let $N_k(r)$ be the number of 
translates of $\mathcal{F}$ contained in
$X_k-U_r(Y_k)$ and $n_k(r)$ the number of 
translates of $\mathcal{F}$ which have a
non-trivial intersection with $U_r(Y_k)$. Set $N_k:=N_k(0)$;
$n_k:=n_k(0)$.
\end{notation}

The next lemma shows that our Definition  \ref{regexhaust} of a regular
exhaustion coincides with the one given by Dodziuk and Mathai \cite{Dodziuk-Mathai(1997)}, with
one exception: we require a lower bound on the injectivity radius of
the boundaries $\boundary X_k$ and control of the covariant derivatives of the
second fundamental form, what they seem to have forgotten (but also
use).

\begin{lemma}\label{longexhaust}
  If $(X_k)_{k \ge 1}$ is a regular exhaustion of $\overline{X}$ as in Definition
  \ref{regexhaust}, then for each $r\ge 0$
 \begin{equation*}
    \lim_{k\to\infty}\frac{\vol(U_r(Y_k))}{\vol(X_k)} = 0.
  \end{equation*}
\end{lemma}
\begin{proof}
  To obtain this we discretize: Choose $\epsilon>0$ such that $4\epsilon$
  is smaller than the injectivity radius, and choose  sets of points
  $P_k\subset Y_k$ such that the balls of radius $\epsilon$ around $x\in P_k$
  are mutually disjoint, but the balls of radius $4\epsilon$ are a covering of
  $Y_k$. Because of bounded geometry (compare the proof of \cite[Lemma
  1.2 in Appendix 1]{Shubin(1992b)}), we find $c_1,c_2>0$ independent
  of $k$ such that
  \begin{equation*}
    c_1 \abs{P_k}\le \area(Y_k) \le c_2 \abs{P_k} .
  \end{equation*}
  The triangle inequality implies $U_r(Y_k)\subset \bigcup_{k}
  B_{r+4\epsilon}(x_k)$. Therefore
  \begin{equation*}
\vol(U_r(Y_k)) \le C_{r+4\epsilon} \abs{P_k}
  \le C_{r+4\epsilon} c_1^{-1} \area(Y_k),
\end{equation*}
  where $C_{r+4\epsilon}$ is a uniform upper
  bound for the volume of balls of radius $r+4\epsilon$ in
  $\overline{X}$ which 
  exists because of bounded geometry. Since we have by assumption
  $\lim_{k \to \infty} \frac{\area(Y_k)}{\vol(X_k)}= 0$, Lemma
  \ref{longexhaust} follows.
\end{proof}

\begin{lemma}\label{boundaryvol}
  If $(X_k)_{k \ge 1}$ is a regular exhaustion of $\overline{X}$ as in Definition
  \ref{regexhaust}, then
  \begin{equation*}
    \lim_{k \to \infty} \frac{\area(\boundary
    X_k\cap \boundary \overline{X})} {\vol(X_k)} ~  = ~
    \frac{\area(\boundary X)}{\vol(X)}.
  \end{equation*}
\end{lemma}
\begin{proof}
Obviously $\vol(\mathcal{F})=\vol(X)$ and
$\area(\mathcal{F}\cap\boundary \overline{X})=\area(\boundary X)$.
Recall that $\mathcal{F}$ is  connected. Suppose that
$\mathcal{F} \cap X_k \not= \emptyset$ and 
$\mathcal{F} \not\subset X_k-U_r(Y_k)$. Then, for each $r$,
$\fundamentalDomain$ must intersect
$U_r(Y_k)$ because otherwise we can find a path in $\mathcal{F}$ connecting a point 
in $X_k$ to a point in $X-X_k$ and this path must meet $Y_k$. Hence we get
for $r \ge 0$
 \begin{align}\label{Nineq1}
   N_k(r)\cdot\vol(X)\le & \vol(X_k) \le (N_k(r)+n_k(r))\cdot\vol(X);
\\
   \label{Nineq2}
   N_k(r)\cdot\area(\boundary X) \le & \area(\boundary X_k\cap\boundary\overline{X})
   \le (N_k(r)+n_k(r))\cdot \area(\boundary X).
 \end{align}
If follows that
\begin{equation}\label{Deep Purple}
 \frac{N_k\cdot \area(\boundary X)}{(N_k + n_k)\cdot\vol(X)} \le
  \frac{\area(\boundary X_k\cap\boundary \overline{X})}{\vol(X_k)}\le
  \frac{(N_k+n_k)\cdot \area(\boundary X)}{N_k \vol(X)}.
\end{equation}
Since $\mathcal{F} \cap U_r(Y_k) \not= \emptyset$ implies 
$\mathcal{F} \subset U_{r + \diam(\mathcal{F})}(Y_k)$, we have
$n_k(r) \cdot \vol(X) \le \vol(U_{r + \diam(\mathcal{F})}(Y_k))$. Therefore
\eqref{Nineq1} implies
\begin{eqnarray*}
\frac{n_k(r)}{n_k(r) + N_k(r)} =
\frac{n_k(r)\vol(X)}{(n_k(r)+N_k(r))\vol(X)}  & \le & 
\frac{\vol(U_{r + \diam(\mathcal{F})}(Y_k))}{\vol(X_k)}.
\end{eqnarray*}
From Lemma \ref{longexhaust}  we conclude
\begin{equation}\label{nineq}
  \lim_{k \to \infty} \frac{n_k(r)}{N_k(r)} =  0.
\end{equation}
Now Lemma \ref{boundaryvol} follows from \eqref{Deep Purple} and
\eqref{nineq}.
\end{proof}

\begin{theorem}\label{randbetticonv}
  If $(X_k)_{k \ge 1}$ is a regular exhaustion of $\overline{X}$ as in Definition
  \ref{regexhaust}, then
  \begin{equation*}
    \lim_{k\to\infty} \frac{b_p(\boundary X_k)}{N_k} =
   \lim_{k\to\infty} \frac{b_p(\boundary X_k)\cdot \vol(X)}{\vol(X_k)} =
    b_p^{(2)}(\boundary \overline{X}).
  \end{equation*}
\end{theorem}
\begin{proof}
  Let $V_k\subset \boundary X_k\cap \boundary \overline{X}$ be the union
  of translates $g\mathcal{F}\cap\boundary\overline{X}$
  for $g \in \Gamma$ such that
  $g\mathcal{F}\subset X_k-Y_k$.
  The number of these translates  $g\mathcal{F}\cap\boundary\overline{X}$
  is just $N_k$. The number $\dot{N}_{m,\delta}$ of ``boundary
  pieces'' appearing in
  \cite{Dodziuk-Mathai(1998)} is bounded by $C_{\delta} \cdot n_k$ for
  a constant $C_{\delta}$ which does not depend on $k$. Because of
  Inequality \eqref{nineq},  $(V_k)_{\ge k}$ is a regular exhaustion of
  $\partial \overline{X}$ in the sense of \cite{Dodziuk-Mathai(1998)} by
  \eqref{nineq}. We conclude from  \cite[Theorem 0.1]{Dodziuk-Mathai(1998)}
  \begin{equation}\label{sconv}
   \lim_{k\to \infty} \frac{b_p(V_k)}{N_k} = b_p^{(2)}(\boundary \overline{X}).
  \end{equation}
  We can thicken $V_k$ inside of $\boundary\overline{X}$ to a regular neighborhood
  $V_k'$. From Proposition \ref{bettibound} we obtain a constant $A$ independent
  of $k$ such that
  \begin{eqnarray}
  b_p(\partial X_k - \interior(V_k'),\partial V_k')  & \le &
  A \cdot \vol(\partial X_k - \interior(V_k')) \nonumber
  \\ & \le &
  A \cdot (\vol(Y_k)+n_k\cdot \vol(\boundary \overline{X}\cap \mathcal{F})).
  \label{Pink Floyd}
  \end{eqnarray}
  We have by excision
  $b_p(\boundary X_k,V_k')= b_p(\boundary X_k-\interior (V_k'),\boundary V_k')$
  and by homotopy invariance $b_p(V_k) = b_p(V_k')$.
  From \eqref{Pink Floyd} and
  the long exact homology sequence of the pair $(\boundary X_k,V_k)$ we conclude
  \begin{equation}
  \label{Roxy music}
  \abs{b_p(\boundary X_k) - b_p(V_k)} \le
  2A \cdot (\vol(Y_k)+n_k\cdot \vol(\boundary \overline{X}\cap \mathcal{F})).
\end{equation}
We get from \eqref{Nineq1} and
\eqref{nineq} (since $\vol(Y_k)/\vol(X_k)\xrightarrow{k\to\infty} 0$
by assumption) that
\begin{equation}\label{eq:Gandalf}
  \lim_{k\to\infty} \frac{2A\left(\vol(Y_k)+n_k \cdot \vol(\boundary
      \overline{X}\cap \mathcal{F})\right)}{N_k} =0.
\end{equation}
  We conclude from \eqref{sconv} and \eqref{Roxy music} and
  \eqref{eq:Gandalf} that 
 \begin{equation}\label{Beatles}
   \lim_{k\to \infty} \frac{b_p(\boundary X_k)}{N_k} = b_p^{(2)}(\boundary \overline{X}).
  \end{equation}
  Now Theorem \ref{randbetticonv} follows from
  \eqref{Nineq1}, \eqref{nineq} and \eqref{Beatles}.

\end{proof}

Remember that the Atiyah-Patodi-Singer index theorem
\cite[Theorem 4.14]{Atiyah-Patodi-Singer(1975a)} and 
\cite[(0.9)]{Cheeger-Gromov(1985a)} and its
$L^2$-version (compare e.g.~\cite{Lueck-Schick(2001b)}) \Kommentar{give 
  precise reference} imply
\begin{equation*}
  \begin{split}
    \frac{\sign(X_k,\boundary X_k)}{\vol(X_k)} =& \frac{1}{\vol(X_k)}\cdot \int_{X_k}
    L(X_k) +
    \frac{\eta(\boundary X_k)}{\vol(X_k)} +
    \frac{1}{\vol(X_k)} \cdot \int_{\boundary X_k} \Pi_L(\boundary X_k),\\
 \frac{\sign^{(2)}(\overline{X},\overline{\boundary X})}{\vol(X)} = &
    \frac{1}{\vol(X)}\cdot\int_{X} L(X)
  +  \frac{\eta^{(2)}(\boundary\overline{X})}{\vol(X)} +
  \frac{1}{\vol(X)}\cdot  \int_{\boundary X} \Pi_L(\boundary X).
\end{split}
\end{equation*}
Here $L(X_k)$ and $L(X)$ denote the Hirzebruch $L$-polynomial, and
$\Pi_L(\boundary X_k)$ and $\Pi_L(\boundary X)$ are local correction terms
which arises because the metric
is not a product near the boundary. We want to show that each of the
individual summands converges for $k\to\infty$ to the corresponding
term for  $\overline{X}$.

The $L$-polynomial is given in terms of the curvature, $\Pi_L$ in
terms of the second fundamental form, therefore both are uniformly
bounded independent of $k$ by some constant $C$.

Moreover, because these are local expressions, the integral over each
translate of the connected fundamental domain $\mathcal{F}$ which is contained
in $X_k-Y_k$ coincides with the corresponding integral on $X$ or
$\boundary X$. Then by splitting the domain of integration appropriately
(as done in the proofs above)
\begin{eqnarray}
\abs{\int_{X_k}L(X_k) - N_k \cdot\int_X L(X)} 
& \le &  
n_k \cdot \vol(X)\cdot C; \label{Iineq1}
\\
\abs{\int_{\boundary X_k} \Pi_L(\boundary X_k) - N_k \cdot 
    \int_{\boundary X}\Pi_L(\boundary X)} 
& \le &
\area(Y_k)\cdot C \nonumber
\\
& & \hspace{2mm} + n_k\cdot\area(\boundary X) \cdot C. \label{Iineq2}
\end{eqnarray}
We conclude from  \eqref{Nineq1} and  \eqref{Iineq1}
\begin{equation}\label{Lest}
  \begin{split}
    & \abs{\frac{1}{\vol(X_k)}\cdot \int_{X_k} L(X_k) -
      \frac{1}{\vol(X)}\cdot \int_X L(X)} \\
    & \le
    \abs{\left(\frac{1}{\vol(X_k)}-\frac{1}{N_k \cdot \vol(X)}\right)\cdot \int_{X_k}
      L(X_k)}\\
    &\qquad + \frac{1}{N_k\cdot \vol(X)}\cdot \abs{\int_{X_k} L(X_k) - N_k\int_X
      L(X)}\\
    & \stackrel{\eqref{Iineq1}}{\le}
    \abs{\frac{1}{\vol(X_k)}- \frac{1}{N_k\vol(X)}}\cdot \vol(X_k)\cdot C +
    \frac{n_k}{N_k}\cdot C \stackrel{\eqref{Nineq1}}{\le} 2C\cdot \frac{n_k}{N_k}.
\end{split}
\end{equation}
We conclude from \eqref{Nineq1}, \eqref{Nineq2}, 
\eqref{Iineq2} and  Lemma \ref{boundaryvol}
\begin{equation}\label{Pest}
  \begin{split}
& \abs{    \frac{1}{\vol(X_k)}\cdot \int_{\boundary X_k}\Pi_L(\boundary X_k) -
  \frac{1}{\vol(X)}\cdot \int_{\boundary X} \Pi_L(\boundary X)}\\
& \le
 \abs{\left(\frac{1}{\vol(X_k)}-\frac{1}{N_k\vol(X)}\right)\cdot \int_{\boundary X_k}
   \Pi_L(\boundary X_k)}\\
 &\hspace{20mm} + \frac{1}{N_k\vol(X)}\cdot \abs{\int_{\boundary
  X_k}\Pi_L(\boundary X_k) - N_k\int_{\boundary X}
  \Pi_L(\boundary X)}\\
 & \stackrel{\eqref{Iineq2}}{\le} \abs{\frac{1}{\vol(X_k)}-
   \frac{1}{N_k\vol(X)}}\cdot \area(\boundary X_k)\cdot C \\
  &\hspace{20mm} +
    \frac{\area{Y_k}}{N_k\vol(X)}\cdot C + \frac{n_k}{N_k}\cdot
  \frac{C\cdot \area(\boundary X)}{\vol(X)} \\
   &\stackrel{\text{\eqref{Nineq1}}}{\le}
   C \cdot \abs{\frac{1}{(n_k + N_k)\cdot \vol(X)} - \frac{1}{N_K\cdot \vol(X)}} \cdot
   \left(\area(Y_k) + \area(\boundary X_k \cap \boundary \overline{X})\right)\\
   & \hspace{20mm} +\frac{\area(Y_k)}{\vol(X_k)}\cdot \frac{C \cdot (N_k + n_k)}{N_k}
  + \frac{n_k}{N_k}\cdot \frac{C\cdot \area(\boundary X)}{\vol(X)}\\
   & \stackrel{\text{\eqref{Nineq2}}}{\le}
    C \cdot \frac{n_k}{N_k \cdot (n_k + N_k)\cdot \vol(X)} \cdot
   \left(\area(Y_k) + (n_k + N_k)\cdot \area(\boundary X)\right)\\
   & \hspace{20mm} + \frac{\area(Y_k)}{\vol(X_k)}\cdot \frac{C \cdot (N_k + n_k)}{N_k}
  + \frac{n_k}{N_k}\cdot \frac{C \cdot \area(\boundary X)}{\vol(X)}\\
   &\stackrel{\text{\eqref{Nineq1}}}{\le} 
   \frac{n_k}{N_k} \cdot \frac{\area(Y_k)}{\vol(X_k)}\cdot \frac{1}{\vol(X)} +  
   \frac{n_k}{N_k}\cdot \frac{\area(\boundary X)}{\vol(X)} \\
   & \hspace{20mm} + \frac{\area(Y_k)}{\vol(X_k)}\cdot \frac{C \cdot (N_k + n_k)}{N_k}
  + \frac{n_k}{N_k}\cdot \frac{C \cdot \area(\boundary X)}{\vol(X)}.
\end{split}
\end{equation}
Since $\lim_{k \to \infty} \frac{\area(Y_k)}{\vol(X_k)} = 0$ by assumption,
from \eqref{nineq}, \eqref{Lest} and \eqref{Pest} follows
\begin{eqnarray*}
\lim_{k \to \infty}  \abs{\frac{1}{\vol(X_k)}\cdot \int_{X_k} L(X_k) -
      \frac{1}{\vol(X)}\cdot \int_X L(X)} & = & 0;
\label{Lconv}
\\
\lim_{k \to \infty} \abs{\frac{1}{\vol(X_k)}\cdot 
\int_{\boundary X_k}\Pi_L(\boundary X_k) -
  \frac{1}{\vol(X)}\int_{\boundary X} \Pi_L(\boundary X)}  & = & 0.
\label{Pconv}
\end{eqnarray*}

It remains to consider the eta-invariants. Because of their
non-local nature this is the most difficult task. The strategy
of the proof of the next
proposition  is similar to the proof of
Remark \ref{etaremark}.

We first recall a few facts about the $\eta$-invariant. Let $D$ be the
tangential signature operator of a $4n-1$-dimensional Riemannian  manifold
$M$, and $\overline M$ a $\Gamma$-covering with lifted signature operator
$\bar D$. Then 
 \begin{equation}   \label{def of eta}
   \eta(M) = \frac{1}{\Gamma(1/2)}
   \int_{0}^\infty t^{-1/2}\tr( D e^{-t D^2})\;dt.
 \end{equation}
and 
\begin{equation} \label{def of eta^{(2)}}
 \eta^{(2)}(\overline{M}) := \frac{1}{\Gamma(1/2)}\int_{0}^\infty t^{-1/2}
 \tr_{\NeumannN\Gamma}(\overline{D} e^{-t\overline{D}^2} )
 \;dt,
\end{equation}
where (with a fundamental domain $\mathcal{F}$ of the covering
$\overline M\to M$)
\begin{eqnarray}
\tr_{\NeumannN\Gamma}(\overline{D} e^{-t\overline{D}^2})
 & = & \int_{\mathcal{F}} \tr_x\left(
 (\overline{D}e^{-t\overline{D}^2})(x,x)\right)dx; \label{gamma_trace and D}
 \\
 \tr(D e^{-t D^2})
 & = &
 \int_{M} \tr_x\left( (D e^{-tD^2})(x,x)\right) dx. \label{trace and D}
\end{eqnarray}

Following Cheeger and Gromov \cite[Section 7]{Cheeger-Gromov(1985a)}
we give an a priori estimate for the large time part of the integrand
defining the $\eta$-invariant. First observe that for $x\ne 0$ we have the
following inequality of functions:
\begin{equation*}\begin{split}
    \abs{\int_T^\infty x t^{-1/2} e^{-tx^2}  \;dt} & = e^{-Tx^2}
    \int_T^\infty \abs{x} t^{-1/2} e^{-x^2(t-T)}\;dt\\
    &\stackrel{t=(u\abs{x}^{-2}+T)} = e^{-Tx^2} \int_0^\infty \abs{x}e^{-u}
    (u\abs{x}^{-2}+T)^{-1/2} \abs{x}^{-2}\;du\\
    &= e^{-Tx^2} \int_0^\infty e^{-u} (u+T\abs{x}^2)^{-1/2} \;du\\
    & \stackrel{T\abs{x}^2\ge 0}{\le}
    e^{-Tx^2}-\chi_{\{0\}} \cdot \sqrt{\pi},
    \end{split}
  \end{equation*}
  where we used $\int_0^\infty u^{-1/2}e^{-u}\;du=\sqrt{\pi}$. For
  $x=0$ obviously $\int_T^\infty xt^{-1/2}d^{-tx^2}\;dt =0$. Hence we get
  for all $x\in\reals$
  \begin{equation*}
    \abs{\int_T^\infty x t^{-1/2} e^{-tx^2}  \;dt} {\le}
    \left(e^{-Tx^2}-\chi_{\{0\}}(x)\right) \cdot \sqrt{\pi},
  \end{equation*}
  where $\chi_{\{0\}}$ is the characteristic function of the set
  $\{0\}$. 
Applying the functional calculus with  $x=\overline D$
 we get
\begin{equation}\label{largeTestoverline{M}}
    \abs{\int_T^\infty t^{-1/2} \tr_{\NeumannN\Gamma}(\overline{D}
     e^{-t\overline D^2})\;dt }
 \le
    \sqrt{\pi} \cdot\tr_{\NeumannN\Gamma}(e^{-T\overline{\Delta}}
     - {\pr}_{\ker\overline{\Delta}}),
\end{equation}
and analogously with $x=D$
\begin{equation}\label{largeTestM_k}
    \abs{\int_T^\infty t^{-1/2} \tr(D
     e^{-t D^2})\;dt }
 \le
    \sqrt{\pi} \cdot \tr(e^{-T\Delta}
     - \pr_{\ker\Delta}).
\end{equation}

\begin{proposition}\label{ametaconv}
  If $(X_k)_{k \ge 1}$ is a regular exhaustion of $\overline{X}$ as in Definition
  \ref{regexhaust} then
  \begin{equation*}
    \lim_{k\to\infty} \frac{\eta(\boundary X_k)}{\vol(X_k)} =
    \frac{\eta^{(2)}(\boundary \overline{X})}{\vol(X)}.
  \end{equation*}
\end{proposition}
\begin{proof}
In the sequel $D[k]$ or $\overline{D}$ is the (tangential)
signature operator
and $\Delta[k]$ or $\overline{\Delta}$ is the differential form Laplacian
on $\boundary X_k$ or $\boundary \overline{X}$, respectively.
Fix $\epsilon > 0$. Choose $T$ such that
  \begin{equation}
    \abs{\tr_{\NeumannN\Gamma}(e^{-T\overline{\Delta}}-
    \pr_{\ker\overline{\Delta}})} \le 
    \frac{\Gamma(1/2)\cdot \vol(X) \cdot \epsilon}{8\sqrt{\pi}}.
    \label{Matthaeus}
  \end{equation}
Put $\boundary X_k^R:=\bigcup_{g\in G\text{ s.t. }U_R(g\mathcal{F})\subset X_k}
(g\mathcal{F}\cap\boundary\overline{X})$.
By Theorem \ref{localize} for the given $T>0$ and $\epsilon>0$ we
find $R>0$ independent of $k$
such that for $0\le t\le T$ and $x\in \boundary X_k^R$.
\begin{eqnarray}\label{kerest1}
\hspace{-7mm} \abs{\tr_x(D[k]e^{-tD[k]^2}(x,x))- \tr_x (
\overline{D} e^{-t\overline{D}^2}(x,x))}
& \le &
\frac{\Gamma(1/2)\cdot \vol(X) \cdot \epsilon}{4\sqrt{T}\cdot \vol(\boundary X)} ;
\\ \label{kerest2}
\hspace{-7mm} \abs{ \tr_x (e^{-t\Delta[k]}(x,x)) -\tr_x
      (e^{-t\overline{\Delta}}(x,x))}
      & \le &
      \frac{\Gamma(1/2)\cdot \vol(X)\cdot \epsilon}
{8\sqrt{\pi}\cdot \vol(\boundary X)} .
\end{eqnarray}
Notice that 
$U_R(g\mathcal{F})\subset X_k \Longleftrightarrow \mathcal{F} \subset X_k-U_R(Y_k)$. 
Hence $\partial X_k^R$ consists of $N_k(R)$ translates of
$\mathcal{F} \cap \boundary \overline{X}$. This implies
$\vol(\boundary X_k^R) = N_k(R) \cdot \vol(\boundary X)$.
From Proposition \ref{bettibound} and \eqref{kerest2} we get for
a constant $A_1$ independent of $k$ (using the fact that
$\overline{\Delta}$ and its kernel are $\Gamma$-equivariant)
\begin{equation}\label{Genesis}
  \begin{split}
    & \abs{ \frac{\tr (e^{-T\Delta[k]})}{N_k(R)}-
   \tr_{\NeumannN\Gamma}(e^{-T\overline{\Delta}})}\\
    & \hspace{5mm} = ~ \abs{\frac{1}{N_k(R)}\cdot \int_{\boundary X_k} 
   \tr_x(e^{-T\Delta[k]}(x,x))\;dx - 
    \int_{\mathcal{F}\cap\boundary \overline{X}}
   \tr_x(e^{-T\overline{\Delta}}(x,x))\;dx}\\
   & \hspace{5mm} \le ~\abs{\frac{1}{N_k(R)}\cdot  \int_{\boundary X_k^R} 
   \tr_x(e^{-T\Delta[k]}(x,x))-
     \tr_x(e^{-T\overline{\Delta}}(x,x))\;dx}\\
   & \hspace{20mm} +
    \abs{\frac{1}{N_k(R)}\cdot \int_{\boundary X_k-\boundary X_k^R} 
    \tr_x(e^{-T\Delta[k]}(x,x))\;dx}\\
   & \hspace{5mm}  = ~
   \frac{\Gamma(1/2)\cdot \vol(X)\cdot \epsilon \cdot \vol(\boundary X_k^R)}
   {N_k(R)\cdot
   8\sqrt{\pi}\cdot \vol(\boundary X)}
   +
   \frac{A_1\cdot \vol(\boundary X_k-\boundary X_k^R)}{N_k(R)}\\
   & \hspace{5mm} \le ~  \frac{\Gamma(1/2)\cdot \vol(X)\cdot \epsilon}{8\sqrt{\pi}} +
   \frac{A_1\cdot \vol(\boundary X_k-\boundary X_k^R)}{N_k(R)}.
   \end{split}
\end{equation}
We conclude from
\eqref{largeTestoverline{M}}, \eqref{largeTestM_k}, \eqref{Matthaeus} and
\eqref{Genesis} (using $\overline D^2=\overline\Delta$ and $D[k]^2=\Delta[k]$)
\begin{eqnarray}
\lefteqn{\left|\frac{1}{N_k(R)}\cdot\frac{1}{\Gamma(1/2)}\cdot 
  \int_{T}^\infty t^{-1/2}\tr(D[k] e^{-t D[k]^2})\;dt
  \right.} & & \nonumber
\\
& & \hspace{8mm} - \left.
\frac{1}{\Gamma(1/2)}\cdot \int_{T}^\infty t^{-1/2}
\tr_{\NeumannN\Gamma}(\overline D e^{-t\overline D^2})
\right| \nonumber
\\
& \le &
\frac{1}{N_k(R)}\cdot \frac{\sqrt{\pi}}{\Gamma(1/2)} \cdot 
\tr(e^{-T\Delta[k]} - \pr_{\ker\Delta[k]}) \nonumber
\\
& & \hspace{8mm} +
\frac{\sqrt{\pi}}{\Gamma(1/2)} \cdot\tr_{\NeumannN\Gamma}(e^{-T\overline{\Delta}}
     - {\pr}_{\ker\overline{\Delta}}) \nonumber
\\
& \le &
\frac{2\sqrt{\pi}}{\Gamma(1/2)} \cdot \tr_{\NeumannN\Gamma}(e^{-T\overline{\Delta}} - 
{\pr}_{\ker \overline{\Delta}}) \nonumber
\\ & & \hspace{8mm}+
\frac{\sqrt{\pi}}{\Gamma(1/2)} \cdot \abs{
\frac{1}{N_k(R)}\cdot\tr(e^{-T\Delta[k]}) -
\tr_{\NeumannN\Gamma}(e^{-T\overline{\Delta}})}\nonumber
\\ & & \hspace{8mm} +
\frac{\sqrt{\pi}}{\Gamma(1/2)} \cdot
\abs{\tr_{\NeumannN\Gamma}({\pr}_{\ker\overline{\Delta}})
-\frac{1}{N_k(R)}\tr(\pr_{\ker\Delta[k]})} \nonumber
\\
& \le & \frac{2\vol(X)\cdot \epsilon}{8} +
\frac{\vol(X)\cdot \epsilon}{8} +
   \frac{\sqrt{\pi}\cdot A_1\cdot \vol(\boundary X_k-\boundary X_k^R) }
{\Gamma(1/2)\cdot N_k(R)}
   \nonumber
\\
& & \hspace{8mm} +
\frac{\sqrt{\pi}}{\Gamma(1/2)}\cdot 
\abs{\tr_{\NeumannN\Gamma}({\pr}_{\ker\overline{\Delta}})
-\frac{1}{N_k(R)} \tr(\pr_{\ker\Delta[k]})}.
\label{Elber}
\end{eqnarray}

From \eqref{gamma_trace and D}, \eqref{trace and D}, \eqref{kerest1}, and
Proposition \ref{smallt} we obtain a constant $A_2$ independent of $k$
such that the following holds:
\begin{eqnarray}
\lefteqn{\left|\frac{1}{N_k(R)}\cdot\frac{1}{\Gamma(1/2)}\cdot 
  \int_0^{T} t^{-1/2}\tr(D[k] e^{-t D[k]^2})\;dt
  \right.} & & \nonumber
\\
& & \hspace{20mm} - \left.
\frac{1}{\Gamma(1/2)}\int_0^{T} t^{-1/2}
\tr_{\NeumannN\Gamma}(\overline{D} e^{-t\overline{D}^2})
\right| \nonumber
\\ & = &
\left|\frac{1}{N_k(R)} \cdot\frac{1}{\Gamma(1/2)}\cdot 
\int_0^T t^{-1/2} \int_{\boundary X_k} \tr_x
    (D[k]e^{-tD[k]^2}(x,x))\;dx\;dt\right.  \nonumber
\\
& & \hspace{20mm} \left.    - \frac{1}{\Gamma(1/2)}\cdot  \int_0^T
  t^{-1/2}\int_{\mathcal{F}\cap\boundary \overline{X}}
  \tr_x(\overline{D}e^{-t\overline{D}^2}(x,x))\;dx\;dt\right| \nonumber
\\
& \le & \frac{1}{N_k(R)} \cdot\frac{1}{\Gamma(1/2)}\cdot \\
& & \hspace{1mm}
 \left|
\int_0^T t^{-1/2} \int_{\boundary X_k^R} \tr_x
    (D[k]e^{-tD[k]^2}(x,x)) -
    \tr_x(\overline{D}e^{-t\overline{D}^2}(x,x))\;dx\;dt\right| \nonumber
\\
& & \hspace{10mm} +
\left| \int_0^T t^{-1/2}
\int_{\boundary X_k - \boundary X_k^R} \tr_x
    (D[k]e^{-tD[k]^2}(x,x))\;dx\;dt\right| \nonumber
\\
& \le & \frac{1}{N_k(R)} \int_0^T t^{-1/2}
\Bigl(\frac{\vol(X) \cdot \epsilon}{4\sqrt{T}\cdot \vol(\boundary X)}  
\cdot \vol(\boundary X_k^R)\\
&& \hspace{17mm}
+ \frac{A_2}{\Gamma(1/2)}
\cdot t^{1/2}\cdot \vol(\boundary X_k - \boundary X_k^R)\Bigr)\; dt.
\nonumber
\\
& \le & \frac{\vol(X) \cdot \epsilon}{8}  +
 \frac{A_2\cdot T}{\Gamma(1/2)}\cdot\frac{\vol(\boundary X_k - \boundary X_k^R)}{N_k(R)}.
\label{Effenberg}
\end{eqnarray}

We conclude from \eqref{def of eta}, \eqref{def of eta^{(2)}}, \eqref{Elber}
and \eqref{Effenberg}
\begin{eqnarray}
\lefteqn{\abs{\eta^{(2)}(\boundary \overline{X}) -
\frac{1}{N_k(R)}\cdot \eta(\boundary X_k)}} & & \nonumber
\\
& \le &
 \frac{\vol(X) \cdot \epsilon}{8} +
 \frac{A_2\cdot T}{\Gamma(1/2)} \cdot\frac{\vol(\boundary X_k- \boundary X_k^R)}{N_k(R)} +
 \frac{3\vol(X)\cdot \epsilon}{8} + \nonumber \\
 & & \hspace{20mm}
\frac{\sqrt{\pi}\cdot A_1\cdot \vol(\boundary X_k-\boundary X_k^R) }{\Gamma(1/2)\cdot N_k}
   \nonumber
\\ & & \hspace{20mm}  +
\frac{\sqrt{\pi}}{\Gamma(1/2)} \cdot
\abs{\tr_{\NeumannN\Gamma}(\pr_{\ker\overline{\Delta}})- \frac{1}{N_k(R)}\cdot 
\tr(\pr_{\ker\Delta[k]})}.
\label{Kahn}
\end{eqnarray}
We get from \eqref{Nineq1}
\begin{eqnarray}
\lefteqn{\abs{\frac{\eta^{(2)}(\boundary \overline{X})}{\vol(X)} -
\frac{\eta(\boundary X_k)}{\vol(X_k)}}} & & \nonumber
\\
& \le &
\frac{N_k(R)}{\vol(X_k)} \cdot
\abs{\eta^{(2)}(\boundary \overline{X}) -
\frac{1}{N_k(R)}{\eta(\boundary X_k)}} +
\abs{\left(\frac{N_k(R)}{\vol(X_k)}- \frac{1}{\vol(X)}\right) \cdot
\eta^{(2)}(\boundary \overline{X})}\nonumber
\\
& \le &
\frac{1}{\vol(X)} \cdot
\abs{\eta^{(2)}(\boundary \overline{X}) -
\frac{1}{N_k(R)}{\eta(\boundary X_k)}} +
\frac{\eta^{(2)}(\boundary \overline{X})}{\vol(X)} \cdot
\frac{n_k(R)}{N_k(R)}.
\label{Babbel}
\end{eqnarray}
We conclude from \eqref{Kahn} and \eqref{Babbel}
\begin{eqnarray}
\lefteqn{\abs{\frac{\eta^{(2)}(\boundary \overline{X})}{\vol(X)} -
\frac{\eta(\boundary X_k)}{\vol(X_k)}}} & & \nonumber
\\
& \le &
 \frac{\epsilon}{8} +
\frac{A_2\cdot T}{\vol(X)\cdot \Gamma(1/2)}
\cdot\frac{\vol(\boundary X_k- \boundary X_k^R)}{N_k(R)} +
\frac{3\epsilon}{8} \nonumber \\
& & \hspace{20mm} +
\frac{\sqrt{\pi}\cdot A_1}{\vol(X)\cdot \Gamma(1/2)}\cdot
\frac{\vol(\boundary X_k-\boundary X_k^R) }{N_k}
   \nonumber
\\ & & \hspace{20mm}  +
\frac{\sqrt{\pi}}{\vol(X)\cdot \Gamma(1/2)} \cdot
\abs{\tr_{\NeumannN\Gamma}(\pr_{\ker\overline{\Delta}})- \frac{1}{N_k(R)}\cdot 
\tr(\pr_{\ker\Delta[k]})}
\nonumber
\\
& & \hspace{20mm}
+ \frac{\eta^{(2)}(\boundary \overline{X})}{\vol(X)} \cdot
\frac{n_k(R)}{N_k(R)}.
\label{Santa Cruz}
\end{eqnarray}
Recall that $\boundary X_k^R$ consists of $N_k(R)$ translates of
$\mathcal{F}\cap\boundary \overline{X}$. The same arguments as above
(using \eqref{nineq}) imply
\begin{eqnarray}
\lim_{k \to \infty} \frac{\vol(\boundary X_k-\boundary X_k^R)}{N_k} & = & 0.
\label{numcount1}
\end{eqnarray}

We get from Theorem \ref{randbetticonv}
  for $k$ sufficiently large
  \begin{equation}
   \label{Dire Strait}
   \lim_{k \to \infty} \frac{\tr(\pr_{\ker\Delta[k]})}{N_k(R)} =
   \lim_{k \to \infty} \frac{b_*(\boundary X_k)}{N_k(R)} =
   b_*^{(2)}(\boundary \overline{X}) =
   \tr_{\NeumannN\Gamma}(\pr_{\ker\overline{\Delta}})
 \end{equation}
 (with the convention $b_*:=\sum_{p\ge 0} b_p$).
From \eqref{nineq}, \eqref{Santa Cruz}, \eqref{numcount1} and
\eqref{Dire Strait} we get the existence of a constant $K(\epsilon)$
such that for all $k \ge K(\epsilon)$
\begin{eqnarray}
\abs{\frac{\eta^{(2)}(\boundary \overline{X})}{\vol(X)} -
\frac{\eta(\boundary X_k)}{\vol(X_k)}}
& \le & \epsilon.
\label{Jancker}
\end{eqnarray}
Now Proposition \ref{ametaconv} follows. This finishes the proof of
Theorem \ref{amconv}.
\end{proof}

\begin{remark}
  Using the symmetry of the tangential signature operator one can
  restrict to $(2n-1)$-forms on the boundary, as explained in
  \cite[Proposition 4.20]{Atiyah-Patodi-Singer(1975b)}. In particular,
  \eqref{Dire Strait} has to hold only for $b_{2n-1}$ (compare also
  \cite{Cheeger-Gromov(1985a),Vaillant(1997)}. 
\end{remark}

\subsection{{Combinatorial version}}
\label{sec:combinatorial}

In this subsection, we prove a combinatorial version of Theorem
\ref{amconv}. It uses the more algebraic techniques employed in Section
\ref{sec:resconv} rather than the heat kernel analysis of Subsection
\ref{sec:analytic}. This way, the result applies to triangulated
rational homology manifolds (with boundary). 

 Let $X$ be a compact
triangulated rational homology manifold with boundary $L$, and of
dimension $4n$. Let
$\overline X$ be a regular covering of $X$ with finitely generated amenable
covering group $\Gamma$. We start by describing the type of exhaustion we
are going to use. Let $\fundamentalDomain$ be
a fundamental domain for the covering $\overline X\to X$,
i.e.~$\fundamentalDomain$
contains exactly one lift of each top-dimensional simplex of $X$. 
For each
simplex $\sigma $ in $X$ choose a lift $\overline \sigma$ in
$\fundamentalDomain$. Let $S$ be a finite system of generators of
$\Gamma$. It gives rise to a left invariant word metric on $\Gamma$. For a
subcomplex $Z\subset \overline X$ and $R>0$, define
\begin{multline*}
  U_R(Z):=\\
  \bigcup_{\sigma \text{simplex of }X}
\{\gamma\overline\sigma\mid \gamma\in\Gamma\text{ and
  }\exists\gamma_1\in\Gamma\text{with
  }d(\gamma-1,\gamma)<R,\gamma_1\overline\sigma\cap Z\ne\emptyset\}. 
\end{multline*}
This depends on the
choice of $S$ as well as the lifts $\overline \sigma$. For each $g\in
\Gamma$, $U_R(g Z)=g U_R(Z)$.

\begin{definition}
  For a simplicial complex $Y$ let $\abs{Y}$ be the total number of
  simplices of $Y$. Similarly, for any subset $W$ of a simplicial
  complex which is a union of open simplices, $\abs{W}$ is the number
  of open simplices in $W$.

  A sequence $X_1\subset X_2\subset \dots \overline X$ of finite
  subcomplexes is called an {\emph{amenable exhaustion}} if
  $\bigcup_{k\in\naturals} X_k=\overline X$ and if for each $R>0$
  \begin{equation*}
    \frac{\abs{U_R(X_k)}}{\abs{X_k}}\xrightarrow{k\to\infty} 1.
  \end{equation*}
  
  It is called a {\emph{balanced}} exhaustion, if for each orbit
  $\Gamma\overline\sigma$ of simplices in $\overline X$
  \begin{equation*}
    \frac{\abs{X_k\cap
        \Gamma\overline\sigma}}{\abs{X_k}}\xrightarrow{k\to\infty}
    \frac{1}{\abs{X}}.
  \end{equation*}

  Denote 
  \begin{equation*}
\tr_k:= \frac{tr_\complexs}{\abs{X_k}}\abs{X};\qquad
\dim_k:=\frac{\dim_\complexs}{\abs{X_k}}\abs{X}.
\end{equation*}
\end{definition}

\begin{theorem}\label{combconv}
  Assume $X$ is a compact simplicial complex triangulating a rational
  homology manifold with boundary the subcomplex $\boundary X$. Assume
  $\overline X$ is a normal covering of $X$ with amenable covering
  group $\Gamma$. Let $X_1\subset X_2\subset \dots$ be subcomplexes forming
  a balanced amenable
  exhaustion of $\overline X$ by rational homology manifolds (with
  boundaries $Y_k$). Then
  \begin{equation*}
    \lim_{k\to\infty} \frac{\sign(X_k,Y_k)}{\abs{X_k}}\abs{X} =
    \sign^{(2)}(\overline X,\boundary\overline X).
  \end{equation*}
\end{theorem}

Before proving this, we investigate the relation between the
Poincar\'e duality maps of one homology manifold being a
codimension-zero subcomplex of
another homology manifold. As an illustration we consider the
following diagram. Let $U\subset M$ be codimension zero
submanifold with boundary $\boundary U$ of a compact manifold $M$. For
the moment assume $\boundary M$ is empty.
\begin{equation*}
  \begin{CD}
    H^p(M) @<<< H^p(M,M-U) @>{\iso}>> H^p(U,\boundary U)\\
    @VV{\cap [M]}V && @VV{\cap [U,\boundary U]}V\\
    H_{n-p}(M) &&&& H_{n-p}(U)\\
    @VV=V && @VVV\\
    H_{n-p}(M) @>>> H_{n-p}(M,M-U) @<{\iso}<<H_{n-p}(U,\boundary U).
  \end{CD}
\end{equation*}
In this diagram, the maps without labels are induced by inclusions and 
the isomorphisms are given by excision. The diagram
commutes because the fundamental class of $M$ is mapped to the
fundamental class of $(U,\boundary U)$ under the composition of the
maps in the lowest row (for $p=0$). 
A corresponding result holds if
$M$ itself has a boundary.

Because we have to apply this in the
$L^2$-setting, we give a chain-level description this diagram. For
this, let $(X,L)$ be an
$n$-dimensional pair of simplicial
complexes  triangulating an oriented rational homology $n$-manifold
$X$ with
boundary $L$. Let $\overline X,\overline L$ be the lifted triangulation of a
normal covering of $X$ ($\overline L$ is the inverse image of $L$ in
$\overline X$) with covering group $\Gamma$. 
Without loss of generality we assume $X$ and $\overline X$ are connected (we
can deal with one component of $X$ at a time, and then the
$L^2$-signature is unchanged if we consider only one component of
$\overline X$).

We first want to describe the (simplicial) $L^2$-chain- and cochain
complexes of $\overline X$.
  Set $\pi:=\pi_1(X)$. 
We have by definition 
\begin{equation*}
  \begin{split}
    C^*_{(2)}(\overline X,\overline L) &= \hom_{\integers\pi }(
    C_*(\tilde X,\tilde L), l^2(\Gamma)),\\
    C^*_{(2)}(\overline X) &=
    \hom_{\integers\pi }( C_*(\tilde X), l^2(\Gamma)),\text{ and}\\
    C_*^{(2)}(\overline X) &= l^2(\Gamma) \tensor_{\integers\pi}
    C_*(\tilde X).
\end{split}
\end{equation*}
Here $\tilde X$ is the
  induced triangulation of the universal covering of $X$, $\tilde L$
  is the inverse image of $L$ in $\tilde X$, and we always
  use the simplicial (co)chain complexes.
 
There are canonical identifications of the simplicial $L^2$-chain and
$L^2$-cochain complexes $C_p^{(2)}(\overline X)$, $C^p_{(2)}(\overline
X)$ with
$L^2$-summable functions on the set of $p$-di\-men\-sional simplices of
$\overline X$, and of $C_p^{(2)}(\overline X,\overline L)$ and
$C^p_{(2)}(\overline X,\overline L)$ with
$L^2$-summable functions on the set of $p$-di\-men\-sio\-nal simplices of
$\overline X$ which do not belong to $\overline L$.

We write $L^2$-functions on the set of simplices of $\overline X$ as
formal sums $\sum \lambda_\sigma\sigma$. Then the identification of
cochains with $L^2$-functions is the anti-linear isomorphism given by $a\mapsto
\sum_{\sigma\in\overline X} \innerprod{1,a(\tilde \sigma)}\sigma$,
where $\tilde \sigma$ is an arbitrary lift of $\sigma$ to the
universal covering. Note that there is a unique projection $p\colon \tilde
X\to \overline X$ since $\overline X$ is a connected normal covering
of $X$. The identification of chains with $L^2$-functions
on the set of simplices is given by $\left(\sum_{g\in
  \Gamma} \lambda_g g\right) \tensor \tilde \sigma\mapsto \sum_{g\in
\Gamma} \lambda_g g \sigma $, where
$\tilde \sigma$ is a simplex in $\tilde X$ (or, for $C_*(\tilde
X,\tilde L)$ of $\tilde X\setminus\tilde L$) and $\sigma=p(\tilde \sigma)$.

Note that this way, in particular we
identify the $L^2$-chain- and cochain groups with each other (via an
anti-linear isomorphism). However, 
this is nothing but the usual isomorphism of a Hilbert
space with its dual. Note that this is not an
isomorphism of chain complexes. Under the
identifications, the chain- and cochain maps induced from the inclusion
of $\overline X$ in $(\overline X,\overline L)$ become the obvious
inclusion and orthogonal projection, respectively.


\begin{lemma}
  Under the above identification, cap-product with the fundamental
  class ---defined on the (co)chain level using the Alexander-Spanier
  diagonal map--- gives a map $C^*_{(2)}(\overline X,\overline L)\to
  C_*^{(2)}(\overline X)$
  which sends
  an $L^2$-function $a$ on the $p$-simplices in $\overline X\setminus \overline
  L$ to
  \begin{equation*}
    \sum_{\overline\sigma\text{ $n$-simplex of $\overline X$}}
    f_{n-p}(\overline\sigma) 
    \innerprod{a,b_p(\overline\sigma)}_{l^2(\overline X)} .
  \end{equation*}
  Here $f_q$, $b_p$ are the front- and back-faces of the corresponding
  dimension, as usual in the Alexander-Spanier diagonal
  approximation. To be able to define this, we choose also a
  $\Gamma$-invariant local ordering of the vertices of $\overline X$,
  e.g.~by lifting such a local ordering from $X$.
\end{lemma}
\begin{proof}
Using the notation introduced above,
$C_*(X,\complexs)$ can be
  identified with $\complexs \tensor_{\integers\pi} C_*(\tilde X)$. For
  each simplex $\sigma$ of $X$ choose a lift $\tilde \sigma$ in
  $\tilde X$.
 Then the fundamental class of $X$ can be written as
  $\sum_{\sigma\in
    X_n} 1\tensor \tilde \sigma$, where $X_p$ denotes
  the collection of $p$-simplices in $X$.

  The Alexander-Spanier cap-product of $a\in \hom_{\integers\pi }(
  C_p(\tilde X,\tilde L), l^2(\Gamma))$ with the fundamental class is then
  given by
  \begin{equation}\label{eq:cap_with_fundamental_class}
     \sum_{\sigma\in
    X_n}  \left(a(b_p(\tilde \sigma))\right)^* \tensor f_{n-p}(\tilde
  \sigma).
\end{equation}
Here $\cdot^*\colon l^2(\Gamma)\to l^2(\Gamma)$ is the standard anti-linear
isomorphism induced from $g\mapsto g^{-1}$ and from complex
conjugation of the coefficients.

Now observe that the function $a=\sum_{\overline\sigma\in\overline
  X}a_{\overline\sigma}\overline\sigma$ is mapped to the
cochain 
\begin{equation*}
\alpha\colon \tilde\sigma\mapsto \sum_{g\in\Gamma}
\overline{a_{g^{-1}p(\tilde\sigma)}} g,
\end{equation*}
and
$\overline{a_{g^{-1}p(\tilde\sigma)}}=\innerprod{g^{-1}p(\tilde\sigma),a}_{l^2}$.

By \eqref{eq:cap_with_fundamental_class}, capping this cochain
$\alpha$ with
the fundamental class gives the chain
\begin{equation*}
  \sum_{\sigma\in X_n} \sum_{g\in\Gamma} \overline{\innerprod{g^{-1}
     p( b_p(\tilde \sigma)),a}} g^{-1}\tensor f_{n-p}(\tilde \sigma).
\end{equation*}
Under our identification, this chain becomes the function 
\begin{equation*}
\sum_{\overline\sigma\in\overline X_n}
\innerprod{a,b_p(\overline\sigma)}f_{n-p}(\overline\sigma),
\end{equation*}
where we use the fact that the family $g^{-1}p(\tilde\sigma)$ for
$g\in\Gamma$ and $\sigma\in X_n$ is exactly the family of all
$n$-simplices of $\overline X$, and the fact that the front-
and back-face maps commute with the action of $\pi$ (and $\Gamma$).
\end{proof}

\begin{lemma}\label{lem:PDfromula}
  Composing the cap-product with the fundamental class with the map
  induced from the inclusion $X\to (X,L)$ 
 we get a
  map $g_{\overline X}:C^*_{(2)}(\overline X,\overline L) \to
  C_*^{(2)}(\overline X,\overline L)$ which under the
  identification with $L^2$-functions on simplices in $\overline
  X\setminus \overline L$ maps such
  a function $a$ (on $p$-simplices) to
  \begin{equation*}
    \sum_{\overline\sigma\text{ $n$-simplex of $\overline X$}}
    f_{n-p}(\overline\sigma) 
    \delta_{\overline L}(f_{n-p}(\overline\sigma))
    \innerprod{a,b_p(\overline\sigma)}_{l^2(\overline X)} .
  \end{equation*}
  Here $\delta_{\overline L}(\overline\sigma)$ is $1$ if $\overline\sigma$ is not
  contained in $\overline L$, and is $0$ if $\overline\sigma\in \overline L$.
\end{lemma}
\begin{proof}
  This is an immediate consequence of the first lemma.
\end{proof}

Assume now $U$ is a compact subcomplex of $\overline X$ which has itself a
subcomplex $V$ (not necessarily contained in $\overline L$) such that
$(U,V)$ triangulates an oriented homology $n$-manifold with boundary
(i.e.~$U$ is a codimension $0$ submanifold with boundary). The above
identifications and formulas apply to the
chain- and cochain complexes of $U$ and $V$ with complex
coefficients. Moreover, observe that these identifications give
canonical embeddings $P_U^*$ of $C^*(U,V,\complexs)$ in
$C^*_{(2)}(\overline X,\overline L)$ and
of $C_*(U,V;\complexs)$ in $C_*^{(2)}(\overline X,\overline L)$. The
corresponding
orthogonal projection is the adjoint $P_U$.

\begin{proposition}\label{restr_eq}
 In this situation, with all the identifications described,
\begin{equation*}
  g_U =    P_U\circ g_{\overline X}\circ P_U^* .
\end{equation*}
\end{proposition}
\begin{proof}
This is implied by the formula of Lemma \ref{lem:PDfromula}. We only
have to make the simple but key
  observation that a top-dimensional simplex of $\overline
  X$ which is not contained in $U$ has no face contained in
  $U\setminus V$ (since in
  the star of an interior point, any two top-dimensional simplices can
  be joined by a sequence of top-dimensional simplices having pairwise
  a common face of codimension $1$. Therefore the star in $\overline
  X$ of an interior point of $U$ can not be bigger than the star in
  $U$).
\end{proof}

Now we are ready to prove Theorem \ref{combconv}. We start with some
auxiliary results we will use. As before, let $\fundamentalDomain$ be
a fundamental domain for the covering $\overline X\to X$. Remember
that for each
simplex $\sigma $ in $X$ we have chosen a lift $\overline \sigma$ in
$\fundamentalDomain$. This way, we get an identification
$C^*_{(2)}(\overline X)=l^2(\overline X)= \oplus_{\sigma\in
  X} l^2(\Gamma) \cdot \overline\sigma$, with subspaces $l^2(X_k)$. Let $P_k^\sigma$ be the orthogonal
projection $l^2(\Gamma)\cdot \overline\sigma \to (l^2(\Gamma)\cdot
\overline\sigma) \cap l^2(X_k)=l^2(X_k\cap \Gamma
\overline\sigma)$. Using the above
identification, $P_k:l^2(\overline X)\to l^2(X_k)$ splits as
$P_k=\diag_{\sigma\in X}(P^\sigma_k)$. For a $\Gamma$-equivariant operator
$ A:C^*(\overline X)\to C^*(\overline X)$ (inducing the operator
$A^{(2)}$ on $C^*_{(2)}$) define $A[k]:= P_k
 A^{(2)} P_k$
(either considered as operator on $l^2(\overline X)$ or on
$l^2(X_k)$). Observe that, if $c:C^*(\overline X)\to
C^*(\overline X)$ 
is the cellular cochain map with adjoint $ c^*$, then $c[k]$ is the
cochain map of $X_k$ with adjoint $c^*[k]$. Note that the
combinatorial Laplacian $\Delta[X_k]=c[k]c^*[k]+c^*[k]c[k]$ in general
does not coincide with $\Delta[k]$ where $\Delta=c c^*+ c^* c$ is the Laplacian of
$\overline X$. By Proposition~\ref{restr_eq} for the Poincar\'e
duality cochain operator we get $g_{X_k}=g[k]$

From this point, the proof of Theorem \ref{combconv} is formally the
same as the proof of Theorem \ref{the: algebraic convergence theorem}
in Section
\ref{sec:resconv}: we have a sequence of operators $g[k]$ and
Laplacians $\Delta[X_k]$ and we have to prove that
\begin{equation*}
  \tr_k\chi_{(a,b)}(\Delta[X_k]
  g[k]\Delta[X_k])\xrightarrow{k\to\infty}\tr_{\NeumannN\Gamma}\chi_{(a,b)}(
  \Delta^{(2)} \overline g\Delta^{(2)})
\end{equation*}
for $(a,b)=(-\infty,0)$ and for $(a,b)=(0,\infty)$. 

We only need the following ingredients, which replace
Lemma \ref{convlemma}, Lemma \ref{balanced1}, Lemma
\ref{lem: normbound}, Lemma \ref{lem:smallEV1}, and Theorem
\ref{betticonv1} in the covering situation, and
the proof given in Section \ref{sec:resconv} goes through.

\begin{lemma}\label{convlemma2}
For $\Gamma$-equivariant linear operators $
h_1,\dots, h_d: \complexs (\overline X)
\to l^2(\overline X)$ (inducing operators $h_k^{(2)}$ on
$l^2(\overline X)$)
and a polynomial $p(x_1,\dots,x_d)$ in non\discretionary{-}{}{-}com\-mut\-ing variables
$x_1,\dots,x_d$ we have
\begin{equation*}
    \tr_{\NeumannN\Gamma}(p(h_1^{(2)},\dots, h_d^{(2)})) =
    \lim_{k\to\infty}tr_k\left(p(h_1[k],\dots, h_d[k])\right) .
\end{equation*}
\end{lemma}
\begin{proof}
  Because of linearity of the traces it suffices
  to study monomials $x_1\dots x_d$. The lemma for $h_1=h_2=\dots h_d$
  and slightly less general
  projections is due to Dodziuk-Mathai
  \cite[Lemma 2.3]{Dodziuk-Mathai(1998)}. An account (with yet another
  slightly different setting) can be found in
  \cite[4.6]{Schick(2001b)}, and
the proof given there carries over with no more
  than obvious changes to the more general situation we are
  considering here.
\end{proof}

\begin{lemma}\label{balanced2}
  For each simplex $\sigma\in X$
  \begin{equation*}
    \tr_k(P_k^\sigma)\xrightarrow{k\to\infty} 1.
  \end{equation*}
\end{lemma}
\begin{proof}
  This is just the definition of a balanced exhaustion.
\end{proof}

\begin{lemma}
  \label{lem: normbound2}
  There is $K\ge 1$ such that for all $k \ge 1$
  \begin{equation*}
    \norm{\Delta^{(2)}}, \norm{\Delta[X_k]},\norm{g^{(2)}},
    \norm{g[k]} \le K.
  \end{equation*}
\end{lemma}
\begin{proof}
  This follows from submultiplicativity of the operator norm and the
  fact that $\norm{P_k}\le 1$ for each $k$, as these are orthogonal
  projections.
\end{proof}

\begin{lemma}\label{lem:smallEV2}
  There is a
constant $C_1>0$ (independent of $k$) such that for
$0 < \epsilon < 1$ and $k \ge 1$
\begin{equation}
\tr_k\left(\chi_{(0,\epsilon]}(\Delta_{2n}[X_k])\right) \le
\frac{C_1}{-\ln(\epsilon)}.
\end{equation}
\end{lemma}
\begin{proof}
  This can either be proved as in \cite[Lemma 2.5]{Dodziuk-Mathai(1998)},
 or we can observe that, since
  $\Delta[X_k]$ is defined over $\integers$, $\ln\Det'\Delta[X_k]\ge
  0$ (compare \cite[Theorem 3.4(1)]{Lueck(1994c)}). Then the inequality follows from Lemma \ref{detcontrol},
  Lemma \ref{lem: normbound2},
  and Lemma \ref{balanced2}.
\end{proof}

\begin{theorem}
  \label{betticonv2}
  The normalized sequence of Betti numbers converges, i.e.~for each $p$
  \begin{equation*}
   \lim_{k\to\infty}    \dim_k(\ker(\Delta_p[X_k])) = \dim_{\NeumannN\Gamma}
   \ker(\Delta^{(2)}_p).
  \end{equation*}
\end{theorem}
\begin{proof}
  This is essentially \cite[Theorem
  0.1]{Dodziuk-Mathai(1998)}. Actually, our exhaustion is slightly
  more general than the ones considered there. But the proof only
  requires Lemma \ref{convlemma2}, Lemma \ref{balanced2}, 
  \ref{lem: normbound2}, and Lemma \ref{lem:smallEV2}; which we have already
  established, and so goes through without changes (compare
  \cite[Section 6]{Schick(2001b)}).
\end{proof}

Now the proof of Theorem \ref{combconv} can be finished as described
above.

\bigskip
Having proved Theorem \ref{combconv}, the question remains open whether in the
 given situation amenable exhaustions by rational homology manifolds
 exist. Amenability is equivalent to the existence F{\o}lner
 exhaustions by subcomplexes without additional structure (used
 e.g.~in \cite{Dodziuk-Mathai(1998)}). We will thicken them to get
 homology manifolds (with boundary), for which the signature is
 defined.
  We thank Steve Ferry who explained to us how to do this
 thickening.

We use the following notation:
\begin{definition}\label{def:star}
  Let $K$ be a simplicial complex with a subcomplex $X$. We define
  $\Star(X)$, the star of $X$, to be the union of the stars of all
  vertices in $X$, where the star of a vertex is the union of all
  closed simplices containing this vertex.

  Denote the barycentric subdivision of $K$ with $K_b$, of $X$ with $X_b$.
\end{definition}

We obviously have:
\begin{lemma}
  In the situation of Definition \ref{def:star},
  $\Star(X)_b=\Star(\Star(X_b))$. 
\end{lemma}

\begin{lemma}
  \label{lemma:thickening}
  Let $(K,L)$ be a triangulated homology manifold (not necessarily
  compact). Let $X'\subset K$ be a subcomplex. Then there exists a
  thickening $X\supset X'$ contained in the star of $X'$, such that $X$ is a
  subcomplex of the barycentric subdivision of $K$ and such that $X$ is
  a rational homology manifold with boundary $Y$. Here $X\cap
  L\subset Y$, but $Y$ is not necessarily contained in $L$.
\end{lemma}
\begin{proof}
  Let $f: K\to\reals$ be a piecewise linear map which is $1$ on $X'$
  and $0$ on the complement of the star of $X'$. Since $f^{-1}(0,1)$
  does not contain a vertex, $1/2$ is a regular value, and therefore
  $X:= f^{-1}([1/2,1])$ will do the job. More specifically,
  $f^{-1}([1/2,1))$ is homeomorphic to the product $f^{-1}(1/2)\times
  [1/2,1)$, since there are no vertices. If $x\in
  X-{f^{-1}(1/2)}$, then it has a neighborhood which is open in $X$ as
  well as in $K$, so it is a manifold point (and will be a boundary
  point whenever it belongs to $L$). All points  $z\in f^{-1}(1/2)$ have
  the neighborhood $U:=f^{-1}([1/2,1))=f^{-1}(1/2)\times [1/2,1)$, and
  no matter how $f^{-1}(1/2)$ looks like, the inclusion $(U-\{y\})\into
  U$ is a homotopy equivalence, so that (by excision) $y$ is a
  boundary point of a rational homology manifold. It remains to
  observe that each path in $X$ is homotopic to a path in
  $X-f^{-1}(1/2)$ such that $X-\boundary X$ does not have more
  connected components than $X$.

  Obviously, we can arrange for $X$ to be a subcomplex of the
  barycentric subdivision of $K$.
\end{proof}

Now we go back to $\overline X$ and construct the
exhaustions we can use. The covering group $\Gamma$ being amenable
means there is a F{\o}lner exhaustion $V_1\subset V_2\subset\dots \Gamma$
with $\bigcup_{k\in \naturals}V_k = \Gamma$
by finite subsets $V_k$, i.e.~$\lim_{k\to\infty}\abs{U_R(\boundary
  V_k)}/\abs{V_k}=0$.
Remember that we have the fundamental domain $\fundamentalDomain$ for
the covering $\overline X\to X$.
If we set $X'_k:= V_k\fundamentalDomain$, then $X'_k$ is an exhaustion
of $\overline X$ by finite subcomplexes as considered in
\cite{Dodziuk-Mathai(1998)}. It is standard that $X_k'$ forms a
balanced amenable exhaustion of $\overline X$. Let $X_k$ be a
thickening of $X_k'$ as
provided by Lemma \ref{lemma:thickening}. Since we want to deal with
(simplicial)
subcomplexes only, we replace $X$ (and $\overline X$) by its
barycentric subdivision. Our main observation is that
$X_k$ is contained in the star of $X'_k$. Fix $R>0$ such that
$\Star(\Star(\fundamentalDomain))\subset U_R(\fundamentalDomain)$. By
the $\Gamma$\discretionary{-}{}{-}in\-vari\-ance of the metric $X_k\subset U_R(X_k')$. Since, on the
other hand, $X_k'\subset X_k$, the sequence $X_k$ forms a balanced
amenable exhaustion of $\overline X$, to which Theorem \ref{combconv}
applies. 


\bibliographystyle{lueck}
\bibliography{lit_A_I,lit_J_Z,signap}

\end{document}